%% file: article.tex
\documentclass[11pt]{article}
\usepackage{amsmath}
\usepackage{amssymb}
\usepackage{amsfonts}
\newenvironment{Proof}{{\begin{minipage}[t]{5.9 in}}{\it Proof.  }}
{{$\spadesuit$}{\end{minipage}}}
\newtheorem{thm}{Theorem}[section]
\newtheorem{lemma}{Lemma}[section]
\newtheorem{prop}{Proposition}[section]
\newtheorem{nott}{Notation}[section]

\newtheorem{conj}{Conjecture}[section]
\newtheorem{cor}{Corollary}[section]

\newtheorem{defn}{Definition}[section]
\newtheorem{remark}{Remark}[section]

\newcommand{\gp}{\mathfrak{p}}
\newcommand{\fr}{\mathop {\rm {Frob_\gp}}\nolimits}
\newcommand{\trace}{\mathop{\rm Trace}\nolimits}
\newcommand{\gal}{\mathop{\rm Gal}\nolimits}
\newcommand{\rank}{\mathop{\rm rank}\nolimits}
\newcommand{\ord}{\mathop{\rm ord}}
\newcommand{\res}{\mathop{\rm res}}
\newcommand{\Q}{{\mathbb{Q}}}
\newcommand{\C}{{\mathbb{C}}}
\newcommand{\Ql}{{\Q}_l}
\newcommand{\Z}{{\mathbb{Z}}}
\newcommand{\p}{{\mathbb{P}}}
\newcommand{\F}{{\mathbb{F}}}
 
\newcommand{\K}{\bar K} 
\newcommand{\bk}{\bar k}

\newcommand{\cS}{\cal S}
\newcommand{\cU}{\cal U}
\newcommand{\cE}{\cal E}
\newcommand{\cT}{\cal T}
\newcommand{\cF}{\cal F}
\newcommand{\cV}{\cal V}
\newcommand{\et}{{\text {\'et}}} 
\newcommand{\tS}{{\tilde {\cal S}}}
\newcommand{\tE}{{\tilde {\cal E}}}
\newcommand{\tF}{{\tilde {\cF}}}
\newcommand{\tT}{{\tilde {\cT}}}
\newcommand{\tD}{\tilde {\Delta}}
\newcommand{\tV}{{\tilde {\cV}}}
\newcommand{\kT}{{\check {\cT}}}
\newcommand{\Fp}{{\F}_\gp}
\newcommand{\bFp}{{\bar {\F}}_\gp}
\newcommand{\qp}{q_\gp}
\newcommand{\Ap}{{A_\gp}(\cE)}
\newcommand{\ap}{a_\gp}
\newcommand{\bp}{b_\gp}
\newcommand{\cp}{c_\gp}
\newcommand{\hK}{\widehat{K}}
\newcommand{\re}{\mathop{\rm Re}\nolimits}
\newcommand{\cOS}{{{\cal O}_{\cS}}}
\newcommand{\cOSS}{{{\cal O}_{{\cS}'}}}
\newcommand{\cOSD}{{{\cal O}_{{\cS},{\Delta}}}}
\newcommand{\tOS}{{{\cal O}_{\tS}}}
\newcommand{\tOSD}{{{\cal O}_{{\tS},{\tD}}}}
\newcommand{\cOE}{{{\cal O}_{\cE}}}

\newcommand{\spec}{\mathop{\rm Spec}\nolimits}
\newcommand{\bspec}{\mathop{\rm {\bf Spec}}\nolimits}
\newcommand{\dvs}{\mathop{\rm Div}\nolimits}
\newcommand{\Hom}{\mathop{\rm Hom}\nolimits}
\newcommand{\pic}{\mathop{\rm Pic}\nolimits}
\newcommand{\alb}{\mathop{\rm Alb}\nolimits}

\newcommand{\pv}{\mathop{\rm Pic^0}\nolimits}
\newcommand{\pvs}{\mathop{\rm Pic^0_{\cS}}\nolimits}
\newcommand{\pvE}{\mathop{\rm Pic^0_E}\nolimits}
\newcommand{\pve}{\mathop{\rm Pic^0_{\cE}}\nolimits}
\newcommand{\pvv}{\mathop{\rm Pic^0_{\cV}}\nolimits}
\newcommand{\pvx}{\mathop{\rm Pic^0_X}\nolimits}
\newcommand{\ns}[1]{{\mathop{\rm NS({\it {#1}})}}}
\newcommand{\mapr}{\smash{\mathop{\longrightarrow}\nolimits}}
\newcommand{\maprlim}[1]{\smash{\mathop{\longrightarrow}\limits^{#1}}}
\newcommand{\maprlow}[1]{\smash{\mathop{\longrightarrow}\limits_{#1}}}
\newcommand{\mapdlim}[1]{\Big\downarrow\rlap
{$\vcenter{\hbox{$\scriptstyle#1$}}$}}

\title{{\bf Arithmetic on Elliptic Threefolds}}
\author{\bigskip \medskip Rania Wazir}

\begin{document}

\maketitle

\bibliographystyle{plain} 
\nocite{*}

\input{artintro.txt}

\input{artbasicdef.txt}

\input{artintsec.txt}

\input{artpicvar.txt}

\input{artsingfib.txt}

\input{artlser.txt}

\clearpage

\appendix

\input{artnot.txt}

\clearpage

\bibliography{artell} 
 
\end{document}

%% file: artintro.txt

\section{Introduction}

\medskip
\par
\noindent
Consider an elliptic curve $E/{\Q}$ given by the Weierstrass equation
\begin{equation}
E: y^2 = x^3 + Ax +B ; \qquad \text {with $A, B \in \Z$},
\label{A1}
\end{equation}
\noindent
and with discriminant locus
$$
\Delta := 4A^3 + 27B^2 \neq 0.
$$

\noindent
The Mordell-Weil Theorem shows that $E({\Q})$, the set of 
rational points on $E$, 
is a finitely-generated Abelian group.  

\medskip
\noindent
If in equation (\ref {A1}) we instead let the coefficients $A, B$ lie in
a polynomial ring over ${\Z}$, then $E$ is no longer defined over ${\Q}$,
but over some function field $K$ of ${\Q}$, and we obtain an 
{\bf elliptic fibration} or {\bf elliptic ${\mathbf n}$-fold}.  
The function-field analogue of the
Mordell-Weil Theorem shows that, in this case also, 
the rational points on $E$ are a 
finitely-generated Abelian group.  The rank of $E(K)$ has been an object of 
intense study and speculation, yet many of its properties, and the 
relation to the underlying geometry of $E$, remain elusive.  The aim of this
paper is to prove a relation between the rank of an elliptic threefold, and
an average of its fibral Frobenius trace values.

\medskip
\noindent 
\begin{thm}
Let $k$ be a number field, ${\gp}$ a prime in $k$, and ${\qp}$ its norm.

\par
\noindent
Let ${\cE} \rightarrow {\cS}$ be a non-split
elliptic threefold defined over $k$.  Then Tate's Conjecture
for ${\cE}/k$ and ${\cS}/k$ implies
$$ \res_{s=1} \sum_{\gp}{-{\Ap}{\frac{\log{\qp}}{{\qp}^s}}}
= \rank{\cE}({\cS}/k).$$
\label{MThm}
\end{thm}

\bigskip

\subsection{History of the problem}

\medskip
\noindent
Consider an elliptic curve $E$ over ${\Q}(T)$.  This is known as an 
{\bf elliptic surface}, and has Weierstrass equation of the form:
$$
E: y^2 = x^3 + A(T)x +B(T) ; \qquad \text {with $A(T), B(T) \in \Z[T]$},
$$
\noindent
with discriminant locus
$$
\Delta(T) := 4A(T)^3 + 27B(T)^2 \not\equiv 0.
$$
\noindent
For each $t \in {\Z}$ and each prime $p$, let
$$
{a_p}(E_t) := 1 + p 
- \#\{\text {rational points in the reduction of $E_t$ mod$p$}\}
$$
\noindent
and let
$$
A_p(E) := {\frac {1}{p}}\sum_{t=1}^{p} {a_p}(E_t)
$$
\noindent
be the average of the ${a_p}$'s over all fibers.

\medskip
\noindent
Then, based on calculations of some non-trivial examples, 
Nagao {\cite {nag}} conjectured that 
$$
\lim_{X \rightarrow \infty} {\frac {1}{X}} 
\sum_{p \leq X} -A_{p}(E)\log{p} = 
\rank E({\Q}(T)).
$$

\par
\noindent
Recently, Rosen and Silverman {\cite {rs}} have derived an analytic version
of Nagao's formula:
$$
\res_{s = 1} \sum_{p} -A_{p}(E){\frac {\log{p}}{{p}^s}} = 
\rank E({\Q}(T)).
$$
\noindent
Assuming Tate's Conjecture, they prove that the analytic version of 
Nagao's formula
holds for elliptic surfaces defined over any number field $k$, and, with a 
mild non-vanishing assumption, that in fact the original Nagao formula holds.

\medskip
\noindent
Theorem {\ref {MThm}} gives an analytic Nagao formula for 
elliptic threefolds.

\subsection{Outline of the Proof}

\medskip
\noindent
The proof proceeds along the following lines:
\begin{enumerate}
\item
We first prove an isomorphism in the cohomology of ${\cE}$ and ${\cS}$:
$$
H^1_{\et}({\cS}/{\bk}, {\Ql}) 
\cong  H^1_{\et}({\cE}/{\bk}, {\Ql})
$$
\noindent
as $\gal({\bk}/k)$-modules.

\item
Next, we need a Shioda-Tate formula for elliptic threefolds.

\item
Find a geometric interpretation of the Frobenius action on the singular fibers
of ${\cE} \mapr {\cS}$.

\item
The number of rational points on a fibered variety 
over ${\Fp}$ can be counted in two ways: we can count the number of rational 
points on the threefold (using the Lefshetz Fixed Point Formula); or we can
count the number of rational points fiber by fiber, and take the sum over
all fibers.  Thus, taking the reduction ${\tE}$ of 
\mbox{${\cE} \mod({\gp})$}, and equating the two expressions for the number 
of rational points on ${\tE}$ will give an equality involving ${\Ap}$.

\item
Re-interpret the equation for ${\Ap}$ in terms of \mbox{$L$-series}, and
apply Tate's Conjecture.

\item
Plug in the results from steps 1, 2, and 3.
\end{enumerate}

\begin{remark}
Steps 1 and 2 are actually proven in the more general case of an
elliptic \mbox{$n$-fold}.  The need for an effective geometric
Tchebotarev theorem requires restriction to the elliptic threefolds 
case for the remaining Steps.
\end{remark}

\begin{remark}
The results in this article formed the bulk of the author's 
Ph.D. thesis {\cite {rw}}.
\end{remark}

%% file: artbasicdef.txt

\section{Basic Definitions and Notation}

\par
\noindent
Most of the results given in this paper, except for the section
on singular fibers, and the concluding section on L-series, 
apply in higher dimensions, and not
just in the case of elliptic threefolds.  Therefore we give here the
general definition of an elliptic n-fold, and where necessary, restrict to 
the case of an elliptic threefold.  This section also contains other basic 
defintions and auxiliary results that will be needed in the rest of 
the article.  For convenience, a list of notation used is included at the end
of the paper.

\medskip
\par
\noindent
\begin{nott}
Let $k$ be a number field, with ring of integers $O_k$, and 
algebraic closure ${\bk}$.
\end{nott}

\subsection{Elliptic N-Folds}

\par
\noindent
\begin{defn}
An {\bf elliptic ${\mathbf n}$-fold defined over ${\mathbf k}$} is a smooth, 
projective variety 
${\cE}/k$ of dimension $n$, together with a proper, flat \mbox{$k$-morphism}
\mbox{$\pi: {\cE} \mapr {\cS}$} to a smooth projective
\mbox{$(n-1)$-dimensional} variety ${\cS}/k$, 
such that the generic fiber is a smooth elliptic curve $E$ defined over
\mbox{$k(\cS)$}, the function field of ${\cS}/k$.
\label{DN0}
\end{defn}

\par
\noindent
In general, {\em morphism} means a morphism defined over
${\bk}$, and {\em \mbox{$k$-morphism}} means a morphism defined over $k$.

\begin{defn}
A {\bf section} \mbox{$\sigma: {\cS} \mapr {\cE}$} is a morphism such
that the composition \mbox{$\pi \circ \sigma = {\rm {id}}_{{\cS}/{\bk}}$}.
\par
\noindent
A {\bf \mbox{${\mathbf k}$-section}} \mbox{$\sigma: {\cS} \mapr {\cE}$} is a 
\mbox{$k$-morphism} such
that the composition \mbox{$\pi \circ \sigma = {\rm {id}}_{{\cS}/k}$}.
\label{DN2}
\end{defn}

\begin{defn}
A {\bf rational section} \mbox{$\sigma: {\cS} \mapr {\cE}$} is a 
rational map (defined over ${\bk}$) such
that the composition \mbox{$\pi \circ \sigma = {\rm {id}}_{{\cU}/{\bk}}$},
for some affine open subset \mbox{${\cU} \subset {\cS}$}.
\par
\noindent
A {\bf \mbox{${\mathbf k}$-rational} section} 
\mbox{$\sigma: {\cS} \mapr {\cE}$} is a 
rational map (defined over $k$) such
that the composition \mbox{$\pi \circ \sigma = {\rm {id}}_{{\cU}/k}$},
for some affine open subset \mbox{${\cU} \subset {\cS}$}.
\label{DN3}
\end{defn}

\par
\noindent
\begin{defn}
An elliptic $n$-fold \mbox{${\cE} \mapr {\cS}$} defined over $k$ is
{\bf split} (over ${\bk}$) if there is an elliptic curve 
${\mathfrak E}$ defined 
over ${\bk}$, and a birational isomorphism (over ${\bk}$) \\
\noindent
\mbox{$\mu:{\cE} \maprlim{\sim} {\mathfrak E} \times_{\bk} {\cS}$}, 
such that the following diagram commutes:
$$
\begin{matrix}
{\cE}   & \maprlim{\mu}  & {\mathfrak E} \times_{\bk} {\cS} \cr
\quad {\scriptstyle \pi} \searrow & & \swarrow {\scriptstyle {\rm proj}_2} \cr
 & {\cS} & \cr
\end{matrix}
$$
\label{DN4}
\end{defn}

\par
\noindent
\begin{nott}
Define ${\mathbf {\cE}/k}$ to be a non-split elliptic $n$-fold  
\mbox{$\pi: {\cE} \mapr {\cS}$} defined over $k$,
with \mbox{$k$-section} \mbox{$\sigma_0: {\cS} \mapr {\cE}$}.  
Assume also that \mbox{${\cE}(k) \neq \emptyset$}.
\label{DN5}
\end{nott}

\medskip
\noindent
\begin{remark}
Note that, because of the section ${\sigma}_0$, 
\mbox{${\cE}(k) \neq \emptyset$} is equivalent to 
\mbox{${\cS}(k) \neq \emptyset$}.  This assumption is made in order
to ensure that ${\pve}$ and ${\pvs}$ are defined over $k$. (See the
discussion in {\cite {sl1}}, pp. 31-33).
\end{remark}

\medskip
\noindent
\begin{nott}
Denote by $(O)$ the image of the section $\sigma_0$ in ${\cE}$, and by
$O$ the corresponding rational point on $E$.
\end{nott}

\begin{defn}
The closed subset
$$
\Delta := \{s \in S \, | \, \text {${\cE}_s$ is not regular}\}
$$
\noindent
is called the {\bf discriminant locus} of ${\cE}$.  $\Delta$ is
a divisor on ${\cS}$.

\par
\noindent
Let $\Delta = {\Delta}_1 + \cdots + {\Delta}_r$ be the
irreducible decomposition of ${\Delta}$.
\label{DN6}
\end{defn}

\par
\noindent
\begin{nott}
$$
\begin{array}{ll}
K           & = k({\cS}) \text {, the function field of ${\cS}/k$.} \\
{\K}        & \text {the algebraic closure of $K$.} \\
{\hK}       & = {\bk}({\cS}) \text {, the function field of ${\cS}/{\bk}$.} \\
{\cOS}      & \text {the structure sheaf of ${\cS}/{\bk}$.} \\
{\cOE}      & \text {the structure sheaf of ${\cE}/{\bk}$.} \\
\end{array}
$$
\label{NT1}
\end{nott}

\par
\noindent
We now make a few observations regarding the elliptic n-fold ${\cE}$.

\begin{prop}
\par
\noindent
\begin{enumerate}
\item[(a)]
\mbox{${\pi}_*({\cOE}) \cong {\cOS}$}.
\medskip

\item[(b)]
The fibers ${\cE}_s$ are connected.
\medskip

\item[(c)]
The requirement that the morphism $\pi$ be flat is equivalent
to requiring that $\dim({\cE}_s)$ be constant for all 
\mbox{$s \in {\cS}$}.

\end{enumerate}
\label{P1}
\end{prop}

\par
\noindent
\begin{Proof}
\begin{enumerate}
\item[(a)]
This follows essentially from the properness of the map $\pi$, and
the existence of a global section $\sigma$.
\par
\noindent
{\bf 1}.  The morphism 
\mbox{$\pi: {\cE} \mapr {\cS}$} can be factored into
\mbox{$\pi = q \circ \pi'$}, where \\
\noindent
\mbox{${\cS}' := {\bspec} \, {\pi}_*{\cOE}$},
\mbox{$\pi': {\cE} \mapr {\cS}'$} is a proper morphism with 
\mbox{${\pi'}_*{\cOE} = {\cOSS}$},
and \mbox{$q: {\cS}' \mapr {\cS}$} is a finite morphism.
(\cite {rh}, proof of Corollary III.11.5).

\par
\noindent

\par
\noindent
{\bf 2}.  By construction, ${\cS}'$ is connected because ${\cE}$ and 
${\cS}$ are.  Furthermore,
the section \mbox{$\sigma: {\cS} \mapr {\cE}$} induces
a section \mbox{${\tilde {\sigma}}: {\cS} \mapr {\cS}'$}, given
by \mbox{${\tilde {\sigma}} := \pi' \circ \sigma$}.  Thus, since
${\tilde {\sigma}}$ maps {\cS} isomorphically onto
its image ${\tilde {\sigma}}(\cS) \subset {\cS}'$,
\mbox{$\dim({\cS}) = \dim({\tilde {\sigma}}(\cS))$}.
But since the morphism $q$ is finite, 
\mbox{$\dim({\cS}) = \dim({\cS}')$ also}.  This implies that 
${\cS}$ is isomorphic to a connected component of ${\cS}'$, and,
since ${\cS}'$ is connected,
the result follows: 
$$
{\cS} \cong {\tilde {\sigma}}(\cS) \cong {\cS}'.
$$

\medskip

\item[(b)]
This follows from part (a), and (\cite {rh}, III.11.3).
\medskip

\item[(c)]
EGA (\cite {ega}, IV.15.4.2).
\end{enumerate} \end{Proof}

\bigskip

\subsection{The ${\hK}/{\bk}$-trace of $E$}

\par
\noindent
In order to state the Mordell-Weil Theorem for function fields, it 
is necessary to recall first the concept of the 
\mbox{${\hK}/{\bk}$-trace} of $E$.

\medskip
\noindent
Let ${\mathfrak k}$ be a field, and let ${\mathfrak F}$ be a 
finitely generated extension of ${\mathfrak k}$, such that ${\mathfrak F}$
is the function field of a variety defined over ${\mathfrak k}$.  Let 
${\mathfrak A}$
be an Abelian variety defined over ${\mathfrak k}$.  An 
\mbox{${\mathfrak F}/{\mathfrak k}$-trace} of ${\mathfrak A}$ is a pair
\mbox{$({\mathfrak t}, {\mathfrak B})$}, consisting of an Abelian variety 
${\mathfrak B}$, and a
homomorphism \mbox{${\mathfrak t}:{\mathfrak B} \mapr {\mathfrak A}$} 
defined over ${\mathfrak F}$, such
that \mbox{$({\mathfrak t}, {\mathfrak B})$} satisfies the following 
universal mapping property:

\medskip
\noindent
If $({\mathfrak r}, {\mathfrak C})$ is another pair with ${\mathfrak C}$ 
defined over ${\mathfrak k}$, and
${\mathfrak r}$ a homomorphism 
\mbox{${\mathfrak r}:{\mathfrak C} \mapr {\mathfrak A}$} defined over 
${\mathfrak F}$,
then there exists a unique homomorphism 
\mbox{${\mathfrak i}:{\mathfrak C} \mapr {\mathfrak B}$} defined over 
${\mathfrak F}$ making the following
diagram commute:

$$
\begin{matrix}
 & {\mathfrak C} \cr
 \quad {\smash{\mathop{\:^{\textstyle\:^{\mathfrak i}}}{\swarrow}}} 
& \mapdlim{\mathfrak r} \cr
 {\mathfrak B} \: \: \maprlow{\mathfrak t} & {\mathfrak A} \cr
\end{matrix}
$$

\bigskip
\noindent
Chow defined and proved the existence of the
\mbox{${\mathfrak F}/{\mathfrak k}$-trace}, and also that the
homomorphism ${\mathfrak t}$ is injective.  See ({\cite {sl1}}, p.138) 
for details.

\medskip
\par
\noindent
Let \mbox{$(\tau, B)$} be the \mbox{${\hK}/{\bk}$-trace} of $E$.
Next, we show that the assumption that the elliptic $n$-fold ${\cE}$
is non-split implies that \mbox{$B = 0$}.

\begin{lemma}
Let ${\cE}$, $E$ be as previously defined, and assume that
the \mbox{${\hK}/{\bk}$-trace} \mbox{$(\tau, B)$} of $E$ is such that 
\mbox{$B \neq 0$}.  Then ${\cE}$ is split over ${\bk}$.
\label{LB}
\end{lemma}

\medskip
\noindent
\begin{Proof}  Since \mbox{$B \neq 0$}, $B$ is an Abelian variety of 
dimension at least one.  Since also the morphism $\tau$ is injective, we 
get
$$
1 \leq \dim B \leq \dim E = 1; 
$$
\noindent
hence \mbox{$\dim B = \dim E = 1$}.
Therefore the morphism \mbox{$\tau: B \mapr E$} is an injective
isogeny, and so
$$
\tau: B \maprlim{\sim} E \quad \text {is an isomorphism over ${\hK}$.}
$$
\noindent
This in turn implies \mbox{${\hK}(B) \cong {\hK}(E)$} as
\mbox{${\hK}$-algebras}.
Since
$$
{\bk}(B \times_{\bk} {\cS}) = {\hK}(B) \cong {\hK}(E) = {\bk}(\cE),
$$
\noindent
this induces a birational isomorphism of varieties
$$
{\cE} \mapr B \times_{\bk} {\cS}
$$
commuting with maps to ${\cS}$; by definition, this means 
${\cE}$ is split.  \end{Proof}

\bigskip

\subsection{$k$-Rational Sections}

\par
\noindent
Let ${\cE}({\cS}/{\bk})$ be the set of rational sections
\mbox{$\rho:{\cS} \mapr {\cE}$}, and 
${\cE}({\cS}/k)$ be the set of \mbox{$k$-rational} sections of ${\cE}$.
${\cE}({\cS}/{\bk})$ and ${\cE}({\cS}/k)$ are groups, and are isomorphic to
$E({\hK})$ and $E(K)$ respectively.  Proofs and further details are exactly
as in the elliptic surfaces case, see ({\cite {js}}, III.3.10).

\par
\noindent
\begin{defn}
${\cE}({\cS}/k)$, the group of \mbox{${\mathbf k}$-rational} 
sections on ${\cE}$, is often called the 
{\bf Mordell-Weil group} of ${\cE}$.
\label{DN7}
\end{defn}

\medskip
\noindent
The next theorem, known as the 
{\bf Mordell-Weil Theorem for function fields},
will show that both ${\cE}({\cS}/k)$ and ${\cE}({\cS}/{\bk})$
have finite rank.  

\begin{thm}[N\'eron - Lang]{\em {\cite {ln}}}
Let ${\mathfrak k}$ be a field, and $\mathfrak K$ 
a function field over $\mathfrak k$.
Let $A$ be an Abelian variety defined over $\mathfrak K$, with
\mbox{${\mathfrak K}/{\mathfrak k}$-trace} $(\mathfrak t, \mathfrak B)$.
Then the group of rational points $A({\mathfrak K})$ modulo the subgroup
${\mathfrak t} {\mathfrak B}({\mathfrak k})$ 
is a finitely generated Abelian group.
\label{TN}
\end{thm}

\medskip
\par
\noindent
Since $B = 0$ by Lemma {\ref {LB}}, this theorem shows that 
$$
E({\hK})/{\tau}B({\bk}) = E({\hK}) \cong {\cE}({\cS}/{\bk})
$$
\noindent
is finitely generated, and therefore  
\mbox{$\rank E({\hK}) = \rank {\cE}({\cS}/{\bk})$} is finite.
Furthermore, ${\cE}({\cS}/k)$ is a subgroup of ${\cE}({\cS}/{\bk})$, so
${\cE}({\cS}/k)$ is also finitely generated;
\mbox{$\rank {\cE}({\cS}/k)$} is often called the {\bf Mordell-Weil rank}
of ${\cE}$.

\medskip
\noindent
\begin{remark}
Actually, ${\cE}({\cS}/k)$ is finitely generated even when 
${\cE}({\cS}/{\bk})$ is not; this follows from Theorem {\ref {TN}}, because
the Mordell-Weil Theorem for number fields shows that $B(k)$ 
is finitely-generated.
\end{remark}

\bigskip

\subsection{Fibral Frobenius Trace Values}

\medskip
\noindent
Having defined one side of the equation in Theorem 1 (the Mordell-Weil rank
of ${\cE}$), we now address the terms on the other side of the 
equation, and make more precise our notion of ``average'' of fibral 
Frobenius trace values:

\begin{nott}
$$
\begin{array}{ll}
{\Fp}       & \text {the residue field of a prime $\gp$ of $O_k$.} \\
{\bFp}  & \text {the algebraic closure of ${\Fp}$.} \\
{\qp}       & \text {the norm of $\gp$, i.e., ${\qp}=\#{\Fp}$.} \\
{\fr}   & \text {the Frobenius morphism over ${\Fp}$.}
\end{array}
$$
\noindent
For a given smooth, projective variety ${\cV}/k$, denote by ${\tV}$ its
reduction mod($\gp$).  
\par
\noindent
${\tV}/{\bFp}$ will be used to denote the variety 
\mbox{${\tV} \times_{\Fp} {\bFp}$}. 
\label{NT2}
\end{nott}

\medskip
\noindent
Some of the traces in the Lefschetz Fixed-Point Theorem 
occur frequently in this
paper, so we make the following definitions:

\medskip
\noindent 
\begin{defn}
For any smooth, projective variety ${\cV}/k$ and any ${\gp}$ 
such that ${\tV}/{\Fp}$ is smooth, let:
$$
\begin{array}{lcl}
{\ap}({\cV})  & := & 
{\trace}({\fr}| H^1_{\et}({\tV}/{\bFp}, {\Ql})),\\
{\bp}({\cV})  & := & 
{\trace}({\fr}| H^2_{\et}({\tV}/{\bFp}, {\Ql})),\\
{\cp}({\cV})  & := & 
{\trace}({\fr}| H^3_{\et}({\tV}/{\bFp}, {\Ql})).\\
\end{array}
$$
\label{DN9}
\end{defn}

\medskip
\noindent
Returning to the case of our elliptic n-fold ${\cE}$, we make the
following definitions:

\medskip
\noindent
\begin{defn} 
For a given point $x \in {\tS}(\Fp)$, let
$$
{\ap}({\cE}_x) := 1 - \#{\tE}_x({\Fp}) + {\qp}m_x,
$$ 
\noindent
where $m_x$ is the number of \mbox{$\Fp$-rational components}
of the fiber ${\tE}_x$.

\medskip
\noindent
These ${\ap}({\cE}_x)$ will be called the {\bf fibral Frobenius trace values}
 of ${\tE}$.
\label{DN10}
\end{defn}

\medskip
\noindent
\begin{remark}
Notice that by the Lefshetz Fixed-Point Theorem, when the fiber ${\tE}_x$
is smooth, definition {\ref {D2}} agrees with definition {\ref {D1}} above.
\label{R1}
\end{remark}

\medskip
\noindent
And finally, the ``average'' of fibral Frobenius trace values:

\medskip
\noindent
\begin{defn} 
$${\Ap} := {\frac{1}{{\qp}^{\scriptstyle (n-1)}}} 
\sum_{x \in {\tS}({\Fp})}{\ap}({\tE}_x).$$
\label{DN11}
\end{defn}

\bigskip

\subsection{Integral Models}

\medskip
\noindent
Consider an embedding of ${\cE}$ into projective space ${\p}^N$.  Since
${\cE}$ is defined over $k$, we can find a finite collection 
$\{ f_1, ..., f_t\}$ of homogeneous polynomials,
with coefficients $\{A_I\}$ in $k$, that define ${\cE}$.  
Singularities on ${\cE}$ are
determined by the simultaneous vanishing of these polynomials, and some
combination of products of their partial derivatives 
(again finitely many, again with
coefficients $\{B_J\}$ in $k$).  Elimination Theory then states that 
there is a set of polynomials $\{ g_h\}_1^u$ in the 
coefficients $\{A_I, B_J\}$, 
such that this is 
equivalent to the simultaneous vanishing of the $g_h$.  Since ${\cE}$
is nonsingular, these $g_h$ do not simultaneously vanish, and hence
there are only finitely many primes ${\gp}$ such that reduction of the
$\{A_I\}$ and $\{B_J\}$ mod(${\gp}$) would give vanishing of the $g_h$. 
Call these the ``bad primes'', and collect them in a set $B$.  The
complement of $B$ is an open set $U \subset O_k$; ${\tE}$, the reduction 
of ${\cE}$ mod(${\gp}$),
is still a smooth, elliptic n-fold (defined over ${\Fp}$)
for all ${\gp} \in U$, and so the model
$$
\begin{array}{c}
{\cE} \\
\downarrow \\
\spec k \\
\end{array}
$$
\noindent
can be extended to an integral model
$$
\begin{array}{c}
{\cE}_U \\
\downarrow \\
U. \\
\end{array}
$$
\noindent
${\cE}/k$ can be thought of as the fiber over the generic point of this
model, while ${\tE}$ is the fiber over the closed point ${\gp} \in U$.

\medskip
\noindent
Similar Elimination Theory arguments can be used to show that the following
hold for all but finitely many primes ${\gp}$:

\begin{enumerate}

\item
${\tS}$ is a smooth $(n - 1)$-fold (defined over ${\Fp}$).
\medskip

\item
${\tilde {\pi}}: {\tE} \mapr {\tS}$ is a proper, flat morphism.

\item
If ${\Delta}'$ is the discriminant locus of 
\mbox{${\tilde {\pi}}: {\tE} \mapr {\tS}$},
then for almost all ${\gp}$, ${\Delta}' = {\tD}$, the reduction
mod(${\gp}$) of ${\Delta}$, and the irreducible decomposition
of ${\Delta}'$ is given by:
$$
{\Delta}' = {\tD}_1 + \dots + {\tD}_r.
$$
\medskip

\item
Let ${\Theta}_{i, j}$ be the irreducible components of 
\mbox{${\pi}^{-1}({\Delta_j})$}, and ${\Theta}_{i, j}'$
be the irreducible components of 
\mbox{${{\tilde {\pi}}}^{-1}({\tD}_j)$}.
Then ${\tilde {\Theta}_{i, j}}$ = ${\Theta}_{i, j}'$.
\medskip

\item
For a fixed $x \in {\cS}$, the number of irreducible components in
${\tE}_x$ is the same as the number of irreducible components in
${\cE}_x$.

\end{enumerate}

\par
\noindent 
By enlarging $B$ if necessary, we will assume that these statements
hold for all ${\gp \in U}$.

%% file: artintsec.txt

\subsection{Intersection Theory}

\medskip
\noindent
One of the main results of the next section will be an isomorphism of the 
Picard Varieties of ${\cE}$ and ${\cS}$.  In the case of elliptic surfaces,
the proof relies on a non-degenerate bilinear pairing 
$$
\dvs({\cE}) \times \dvs({\cE}) \mapr {\Z}
$$
given by Intersection Theory on surfaces.  (Details can be found in
{\cite {js}} Proposition III.8.2).  In the case of higher-dimensional
varieties, it is no longer possible to get a pairing into ${\Z}$.  However,
we will show that it suffices to have a pairing with a notion of 
``positivity.''  This pairing will be the subject of this
section.  For a review of the necessary results
on Intersection Theory of higher-dimensional varieties, see
\mbox{\cite{rh}, Appendix A} or \mbox{{\cite{wf}}}.  For ease of reference,
we list here only the most necessary results.

\medskip
\noindent
For the rest of this section, assume all varieties are non-singular, 
projective over an algebraically closed field ${\kappa}$.

\subsubsection{Summary of Basic Results}

\par
\noindent
\begin{defn}
A cycle \mbox{$\Gamma = \sum a_i[V_i]$} is {\bf non-negative} if
\mbox{$a_i \geq 0$} for all $i$; in that case, we write 
\mbox{$\Gamma \geq 0$}.  
$\Gamma$ is {\bf positive}, written \mbox{$\Gamma > 0$}, 
if in addition, \mbox{$a_i > 0$} for some $i$.
\label{D1}
\end{defn}

\begin{defn}
$\Gamma$ is {\bf non-positive}, written \mbox{$\Gamma \leq 0$}, 
if \mbox{$a_i \leq 0$} for all i.
$\Gamma$ is {\bf negative}, written \mbox{$\Gamma < 0$}, 
if in addition, \mbox{$a_i < 0$} for some $i$.
\label{D2}
\end{defn}

\medskip
\noindent
\begin{defn}
A cycle class \mbox{$\alpha \in A_r(X)$} is {\bf non-negative}, written 
\mbox{$\alpha \geq 0$}, if there is a non-negative $r$-cycle ${\Gamma}$
such that \mbox{$\Gamma \geq 0$} and \mbox{cl$(\Gamma) = \alpha$}.  

\par
\noindent
Similarly for $\alpha$ {\bf positive} (\mbox{$\alpha > 0$}),
$\alpha$ {\bf non-positive} (\mbox{$\alpha \leq 0$}), and 
$\alpha$ {\bf negative} (\mbox{$\alpha < 0$}).
\end{defn}

\medskip
\noindent
\begin{thm}
If \mbox{$f: X \mapr W$} is a proper morphism, and ${\Gamma}$ is an $r$-cycle
on $W$ which is algebraically equivalent to zero, then 
$f_*{\Gamma}$ is algebraically equivalent to zero on $X$.
\label{I1}
\end{thm}

\par
\noindent
\begin{Proof} See Fulton(\cite{wf}, Proposition 10.3). \end{Proof} 

\medskip
\noindent
\begin{thm}
If \mbox{$f: X \mapr W$} is a morphism, and ${\Gamma}$ is an $r$-cycle
on $W$ which is algebraically equivalent to zero, then 
$f^*{\Gamma}$ is algebraically equivalent to zero on $X$.
\label{I2}
\end{thm}

\par
\noindent
\begin{Proof} See Fulton(\cite{wf}, Example 19.3.9, p.390). \end{Proof}

\bigskip

\subsubsection{A Nondegenerate Pairing}

\begin{defn}
Define a symmetric, bilinear pairing 
$$
\langle \cdot, \cdot \rangle:{\pic}({\cE}) \times {\pic}({\cE}) 
\mapr {\pic}({\cS})
$$
\noindent
via 
$$
\langle {\Lambda}, {\Upsilon} \rangle := {\pi}_*({\Lambda}.{\Upsilon})
$$
\noindent
for any \mbox{${\Lambda}, {\Upsilon} \in {\pic}({\cE})$}.

\medskip
\par
\noindent
If $C, D$ are divisors in ${\dvs}({\cE})$, let 
\mbox{$\langle C, D \rangle := \langle cl(C), cl(D) \rangle$}.
\label{D11}
\end{defn}

\medskip
\noindent
The following Proposition is very similar to the Elliptic Surfaces case.  
However, there are enough differences due to the new pairing, that we 
give a proof below.

\begin{prop}
Let \mbox{$D \in {\dvs}({\cE})$} be a fibral divisor, and 
\mbox{$G \in {\dvs}({\cS})$}.
Then
\begin{itemize}
\item[(a)]
$\langle D, {{\pi}^*}(G) \rangle = 0$
\item[(b)]
$\langle D, D \rangle \leq 0$
\item[(c)]
If $\langle D, D \rangle = 0$ then 
$D \in {\pi}^*({\dvs}({\cS}))$.
\end{itemize} 
\label{I3}
\end{prop}

\medskip
\noindent
\begin{Proof}
\begin{itemize}
\item[(a)]
This follows from the Projection Formula, 
 once we note that for any fibral divisor $D$, \mbox{${\pi}_*(D) = 0$}.
\begin{eqnarray*}
\langle D, {{\pi}^*}(G) \rangle & = &
{\pi}_*(cl(D).cl({{\pi}^*}{G})) \\
  & = & {\pi}_*(D.{{\pi}^*}(G)) \\
  & = & {\pi}_*(D).G  \\
  & = & 0.G \\ 
  & = & 0.
\end{eqnarray*}

\medskip

\item[(b)]
Write 
$$D = D_1 + D_2 + ... + D_n,$$
\noindent
where each $D_i$ is contained in a 
different fiber.  Then \mbox{$\langle D_i, D_j \rangle = 0$} 
for $i \neq j$, which 
implies that 
$$\langle D, D \rangle = \langle D_1, D_1 \rangle 
+ \langle D_2, D_2 \rangle + ... + \langle D_n, D_n \rangle.$$
\noindent
Thus, it suffices to prove the theorem for each $D_i$, and we
can assume that \mbox{$D \subset {\pi}^*G$}, for some 
irreducible \mbox{$G \in Div(S)$}.  

Let 
$$
F := {\pi}^*G = \sum_{i=0}^{t} n_i{\Gamma}_i
$$
\noindent
be the irreducible decomposition of $F$.  Notice that
\mbox{$n_i \geq 0$} for all $i$.  Furthermore, since 
\mbox{$D \subset F$}, $D$ can be written as 
$$
D = \sum_{i=0}^{t} a_i{\Gamma}_i.
$$

\noindent
Rewrite $D$ as
$$
D = \sum_{i=0}^t ({\frac{a_i}{n_i}}){n_i}{\Gamma}_i,
$$
\noindent
and define another fibral divisor $D'$ by
$$
D' = \sum_{i=0}^t ({\frac{a_i}{n_i}})^2{n_i}{\Gamma}_i.
$$

\noindent
By part (a), we know that 
\mbox{$\langle D', F \rangle = 0$} and thus:
\begin{eqnarray*}
-2{\langle D, D \rangle} & = & {\langle D', F \rangle} -2{\langle D, D \rangle}
 + {\langle F, D' \rangle} \\
  & = & \sum_{i,j = 0}^t {{\frac{a_i^2}{n_i^2}}
{\langle n_i{\Gamma}_i, n_j{\Gamma}_j \rangle}}
-2\sum_{i,j = 0}^t {{\frac{a_ia_j}{n_in_j}}
{\langle n_i{\Gamma}_i, n_j{\Gamma}_j \rangle}} \\
  &   & + \sum_{i,j = 0}^t {{\frac{a_j^2}{n_j^2}}
{\langle n_i{\Gamma}_i, n_j{\Gamma}_j \rangle}} \\
  & = & \sum_{i,j = 0}^t {(\frac{a_i}{n_i} - \frac{a_j}{n_j})^2
{\langle n_i{\Gamma}_i, n_j{\Gamma}_j \rangle}} \\
\end{eqnarray*}

\noindent
The terms in the above equation are zero when $i=j$; it therefore reduces to:
\begin{equation}
{\langle D, D \rangle} = -{\frac{1}{2}} \sum_{{{i,j = 0} \atop {i \neq j}}}^t 
{(\frac{a_i}{n_i} - \frac{a_j}{n_j})^2
{\langle n_i{\Gamma}_i, n_j{\Gamma}_j \rangle}}
\label{B1}
\end{equation}

\noindent
Now for $i \neq j$, ${\Gamma}_i$ and ${\Gamma_j}$ are distinct irreducible 
divisors, and so 
$${\langle {\Gamma}_i, {\Gamma}_j \rangle} \geq 0 \quad 
\text {for all $i \neq j$.}$$
\noindent
Since, as previously mentioned, \mbox{$n_i, n_j \geq 0$}, this implies that
$${\langle n_i{\Gamma}_i, n_j{\Gamma}_j \rangle} \geq 0 \quad 
\text {for all $i \neq j$,}$$
\noindent
and therefore
\mbox{${\langle D, D \rangle} \leq 0$}, proving part (b).
\medskip

\item[(c)]
To prove part (c), now assume that \mbox{${\langle D, D \rangle} = 0$}. 
Plugging into equation (\ref{B1}), this gives 
$$
\frac{a_i}{n_i} = \frac{a_j}{n_j} \quad \text {for all $i, j$ such that 
$\langle {\Gamma}_i, {\Gamma}_j \rangle > 0$.}
$$

\noindent
However, the fibers of $\pi$ are connected by Proposition {\ref {P1}}.b.
Thus,  given any two irreducible 
components ${\Gamma}_i$ and ${\Gamma_j}$, there is a sequence of components
$$
{\Gamma}_i = {\Gamma}_{i_0}, {\Gamma}_{i_1}, ..., {\Gamma}_{i_r} = {\Gamma}_j
$$
\noindent
such that ${\Gamma}_{i_k}$ has non-empty intersection with 
${\Gamma}_{i_{k+1}}$.
Hence 
$$\frac{a_i}{n_i} = \frac{a_j}{n_j} \quad \text {for all $i, j$.}$$
\noindent
But the section \mbox{$\sigma: S \mapr {\cE}$} intersects with
exactly one irreducible component of ${\pi}^*G$, say $\Gamma_0$, which
therefore has multiplicity one:
i.e., \mbox{$n_0 = 1$}.
This implies
$$\frac{a_i}{n_i} = a_0 \in {\Z} \quad \text {for all $i$}.$$
\noindent
Plugging this into the irreducible decomposition of $D$ gives:
\begin{eqnarray*}
D & = & \sum_{i=0}^t {a_i{\Gamma_i}} \\
  & = & \sum_{i=0}^t {{\frac{a_i}{n_i}}n_i{\Gamma}_i}\\
  & = & \sum_{i=0}^t {a_0n_i{\Gamma}_i}\\
  & = & a_0\sum_{i=0}^t {n_i{\Gamma}_i}\\
  & = & a_0F \in {\pi}^*({\dvs}({\cS})),
\end{eqnarray*}

\noindent
This completes the proof of part (c). \end{itemize} \end{Proof}

\bigskip

\subsection{Picard Varieties}

\medskip
\noindent
We will adopt the following notation for varieties $X/{\kappa}$, ${\kappa}$
an algebraically closed field:
\begin{nott} 
$$
\begin{array}{lcl}
{\pic}_{X} &  & \text {the Picard Scheme of $X/{\kappa}$.} \\
{\pic}(X) & = & {{\pic}_X}({\bar {\kappa}}) 
\text {, the Picard Group of $X/{\bar {\kappa}}$.}\\
{\pic}(X/{\kappa}) & = & {{\pic}_X}(\kappa) 
\text {, the Picard Group of $X/k$.}\\
{\pvx} &  & \text {the Picard Variety of $X/{\kappa}$.} \\
{\pv}(X) & = & {{\pvx}}({\bar {\kappa}}). \\
{\pv}(X/{\kappa}) & = & {{\pvx}}({\kappa}). \\  
{\ns{X}} & = & {\pic}(X)/{\pv}(X) \text {, the Neron-Severi Group of 
$X/{\bar {\kappa}}$.} \\
{\ns{X/{\kappa}}} & = & {\pic}(X/{\kappa})/{\pv}(X/{\kappa}) 
\text {, the Neron-Severi Group of $X/{\kappa}$.} 
\end{array}
$$
\label{D4}
\end{nott}

\medskip
\noindent
\begin{defn}
The trivial part of ${\ns{\cE}} \otimes \Q$, denoted ${\cT}$, is the 
subspace generated by the image of the zero section $(O)$, and by all 
geometrically irreducible components of the fibral divisors.  
\par
\noindent
Denote by ${\cF}$ the subspace of ${\cT}$ generated by the
non-identity components of the fibral divisors, where the identity
component of a fibral divisor is the component intersecting $(0)$.  
\label{D5}
\end{defn}

\par
\noindent
Note that ${\cT}$ is generated by $(0)$, ${\pi^*}{\ns{\cS}}$, and ${\cF}$.

\par
\noindent
\begin{remark}
For all but finitely many primes $\gp$, the trivial part of
\mbox{$\ns{{\tE}/{\bFp}} \otimes \, \Q$} is isomorphic to ${\tT}$, the
reduction of ${\cT} \mod(\gp)$.  This follows from Section 2.5,
statement 3, and hence holds for all primes ${\gp} \in U$.  
\label{R7}
\end{remark}

\bigskip

%% file: artpicvar.txt

\section{An Isomorphism in Cohomology}

\medskip
\noindent
Our main goal in this section will be to prove
\medskip
\noindent
\begin{thm}
$$
H^1_{\et}({\cS}/{\bk}, {\Ql}) 
\cong  H^1_{\et}({\cE}/{\bk}, {\Ql}).
$$
\noindent
as $\gal({\bk}/k)$-modules.
\label{T8}
\end{thm}

\noindent
This can be proven by showing first that \mbox{${\pve} \cong {\pvs}$}, and
combining this with the fact that, for any variety ${\cV}/k$, 
\mbox{$H^1_{\et}({\cV}/{\bk}, {\Ql}) \cong H^1_{\et}({\pvv}/{\bk}, {\Ql})$}.

\medskip
\noindent
The isomorphism of Picard Varieties will follow quite easily from the 
following theorem: 

\begin{thm}
There is an exact sequence of Abelian Groups:
\begin{equation}
0 \mapr {\pv}({\cS}) \mapr {\pv}({\cE}) \mapr B({\bk}).
\label{B6}
\end{equation}
\label{T2}
\end{thm}

\par
\noindent
\begin{Proof}  This proof is based on Raynaud's proof of the
elliptic surface case (see {\cite {sh1}}, Theorem 2).

\medskip
\noindent
\begin{enumerate}
\item
\par
\noindent
The morphism \mbox{${\pi}: {\cE} \rightarrow {\cS}$} induces a map 
$$
{\pi}^*: {\pic}_{\cS} \rightarrow {\pic}_{\cE}.
$$  
\noindent
Since ${\pi}^*$ preserves algebraic equivalence (Theorem \ref {I2}), we 
get a morphism 
$$
{\alpha} := {\pi}^*|_{\pvs}: {\pvs} \rightarrow {\pve}.
$$
\noindent
Furthermore, restriction to the generic fiber induces a morphism 
\begin{equation}
{\psi}: {{\pve} {\times}_{\cS} {\cS}} \rightarrow {\pvE} \cong E.
\label{B9}
\end{equation}

\par
\noindent
But by the universal mapping property of the 
\mbox{${\hK}/{\bk}-$trace $(\tau,B)$} of $E$, the map $\psi$
factors through $B$; i.e., there is a unique homomorphism
\mbox{$\beta: {\pve} \mapr B$} such that
$\psi = \tau \circ \beta$.
Thus, we have a sequence of morphisms

\begin{equation}
{\pvs} \maprlim{\alpha} {\pve} \maprlim{\beta} B
\label{B2}
\end{equation}

\noindent
and it remains to show that, as maps on the ${\bk}$-points, this is a 
short exact sequence of Abelian groups:

\begin{equation}
0 \mapr {\pv}({\cS}) \maprlim{\alpha} {\pv}({\cE}) \maprlim{\beta} B({\bk}).
\label{B3}
\end{equation}

\medskip

\item
\par
\noindent
The injectivity of ${\alpha}$ follows from the existence of the global section
\mbox{$\sigma: {\cS} \mapr {\cE}$};
\par
\noindent
Since by definition \mbox{${\pi} \circ {\sigma} = id_{\cS}$}, 
this gives
$$id_{{\pic}({\cS})} = ({\pi} \circ {\sigma})^* = {\sigma}^* \circ {\pi}^*;$$
\noindent
therefore ${\pi}^*$, and hence also ${\alpha}$, is injective.

\medskip

\item
\noindent
To show exactness at 
the middle, note first that \mbox{${\psi} \circ {\alpha} = 0$}; 
this is true because restriction
to the generic fiber sends fibral divisors in ${\pv}({\cE})$ to zero, 
and $\alpha$ sends 
${\pv}({\cS})$ to fibral divisors in ${\pv}({\cE})$.  

\par
\noindent
But if \mbox{${\psi} \circ {\alpha} = 0$}, then also
\mbox{${\beta} \circ {\alpha} = 0$} because $\tau$ is injective.
\par
\noindent
This shows \mbox{Im($\alpha$) $\subset$ Ker($\beta$)}.

\medskip

\item
\noindent
Finally, to show \mbox{Ker($\beta$) = Im($\alpha$)},
take any \mbox{$0 \neq \Gamma \in$ Ker($\beta$)}.  
Then \mbox{$\Gamma = cl(D)$} for some
divisor $D$ with \mbox{$D|_E = 0$}, 
and thus $D$ must be a fibral divisor.  
Since \mbox{$cl(D) \in {\pv}({\cE})$}, it follows that 
\mbox{$\langle D, F \rangle = 0$} for every divisor $F$ on ${\cE}$,
and in particular, \mbox{$\langle D, D \rangle = 0$}.
By Theorem (\ref{I1}.c), this implies
\mbox{$D \in {\pi}^*({\dvs}({\cS}))$}.

\medskip
\noindent
Thus, we have $D \sim_{alg} 0$, and $D = {\pi}^*C$ for some
$C \in \dvs({\cS}).$  It remains to show that $C \sim_{alg} 0$.
But by Theorem {\ref {I2}}, $\sigma^*$ preserves algebraic equivalence.
Therefore, $\sigma^*D \sim_{alg} 0$.
Since
\begin{eqnarray*}
\sigma^*D & = & \sigma^*\pi^*C \\
 & = & (\pi \circ \sigma)^*C \\
 & = & (id_{\pic({\cS})})C \\
 & = & C,
\end{eqnarray*}
\noindent
this shows that $C \sim_{alg} 0$, and therefore 
\mbox{${\Gamma} = cl(D) \in {\pi}^*({\pv}({\cS}))$}, proving that 
\mbox{Ker($\beta$) = Im($\alpha$)}.  

\end{enumerate}
\end{Proof}

\noindent
Now we are ready to prove the isomorphism of Picard Varieties:

\begin{thm}
${\pvs}$ and ${\pve}$ are isomorphic as Abelian Varieties over $k$.
\label{T3}
\end{thm}

\medskip
\noindent
\begin{Proof}  
\begin{enumerate}
\item
Since the elliptic threefold ${\cE}$ is non-split, the 
\mbox{${\hK}/{\bk}$-trace $B$} must be trivial (Lemma {\ref {LB}}).  
From the exact sequence 
(\ref{B6}), it then
follows that ${\pv({\cS})}$ and ${\pv({\cE})}$ are isomorphic as groups.  
\medskip

\item
If \mbox{${\pv({\cS})} \cong {\pv({\cE})}$}, then 
\mbox{$\pvs \cong \pve$}.

\par
\noindent
For any abelian variety $A/k$, and any integer $l$, let
$$
A_l := \{ a \in A \, | \, la = 0 \}.
$$
\noindent
Then by {\cite {dm1}}, Proposition II.6, p.64, 
$$
A_l \cong {{\Z}/{l{\Z}}}^{2d_A},
$$
\noindent
where $d_A := \dim(A)$.

\par
\noindent
Therefore, \mbox{${\pv({\cS})} \cong {\pv({\cE})}$} implies
\mbox{${\pv({\cS})}_l \cong {\pv({\cE})}_l$}, and so
$$
\dim(\pv({\cS})) = \dim(\pv({\cE})).
$$

\par
\noindent
Since $\alpha$ is a map between two Abelian varieties of the same dimension,
and with finite (in fact, trivial) kernel, this implies that $\alpha$ is an 
isogeny. ({\cite {jm2}}, Proposition 8.1.c)

\par
\noindent
But then 
$$
\pve \cong \pvs/Ker(\alpha), 
\quad \text {by {\cite {dm1}}, Corollary III.10.1, p.118,}
$$
\noindent
and since $\ker(\alpha) = 0$, this gives
$\pve \cong \pvs$.
\medskip

\item
Since $\pvs$, $\pve$, and the map $\alpha$ are all defined over $k$,
the isomorphism is in fact an isomorphism over $k$.

\end{enumerate}
\end{Proof}

\medskip
\noindent
\begin{prop}
Let ${\cV}$ be a smooth, projective variety defined over $k$. Then
$$
H^1_{\et}({\cV}/{\bk}, {\Ql}) \cong H^1_{\et}({\pvv}/{\bk}, {\Ql}).
$$
\medskip
\label{PVH}
\end{prop}

\newcommand{\mt}{\mathfrak T}

\medskip
\noindent
\begin{Proof}
This seems to be a well-known result, but since we could not find a 
reference, we give a proof below.

\medskip
\noindent
Denote by ${\alb}_{\cV}$ the Albanese variety of ${\cV}$, and let
$$
\mu:{\cV} \mapr {\alb}_{\cV}
$$
be the Albanese map of ${\cV}$. 
This map is defined over $k$ if ${\cV}(k)$ is not empty - 
({\cite {sl2}}, pp.31-32).
Then $\mu$ induces an isomorphism in cohomology ({\cite {gh}}, p.331)
$$
H_1({\cV}, \Z)/{\mathfrak T} \maprlim{\sim} 
H_1({\alb}_{\cV}, \Z),
$$
\noindent
where ${\mathfrak T}$ is the torsion subgroup of $H_1({\cV}, \Z)$.
Rewrite this as a short exact sequence
$$
0 \mapr {\mt} \mapr H_1({\cV}, \Z) \mapr H_1({\alb}_{\cV}, \Z) \mapr 0.
$$
Taking $\Hom(\: \cdot \:, {\Z}/{l^r{\Z}})$ of the sequence gives:
$$
0 \mapr \Hom(H_1({\alb}_{\cV}, \Z), {\Z}/{l^r{\Z}}) 
\mapr \Hom(H_1({\cV}, \Z), {\Z}/{l^r{\Z}}) \mapr \Hom({\mt}, {\Z}/{l^r{\Z}}).
$$
Since 
\mbox{$\Hom(H_1(X, \Z), G) \cong H^1(X, G)$} for any finite Abelian group $G$ 
and any smooth variety $X$, we get
$$
0 \mapr H^1({\alb}({\cV}), {\Z}/{l^r{\Z}})
\mapr H^1({\cV}, {\Z}/{l^r{\Z}}) \mapr \Hom({\mt}, {\Z}/{l^r{\Z}}).
$$
\noindent
Furthermore, by the \`etale cohomology compatibility theorems 
({\cite {jm}}, Theorem III.3.12, p.117), we have
\mbox{$H^1(X, G) \cong H^1_{\et}(X, G)$} for any finite Abelian group $G$ and
any smooth variety $X/k$:
$$
0 \mapr H^1_{\et}({\alb}_{\cV}, {\Z}/{l^r{\Z}})
\mapr H^1_{\et}({\cV}, {\Z}/{l^r{\Z}}) \mapr \Hom({\mt}, {\Z}/{l^r{\Z}}).
$$
\noindent
Taking inverse limits as $r$ goes to infinity, and noting that inverse limits
of exact sequences of finite groups are exact ({\cite {jm}}, p.165):
$$
0 \mapr H^1_{\et}({\alb}_{\cV}, {\Z}_l) \mapr H^1_{\et}({\cV}, {\Z}_l) 
\mapr \lim_{\leftarrow} \Hom({\mt}, {\Z}/{l^r{\Z}}),
$$
\noindent
where by definition, 
$$
H^1_{\et}({\cV}, {\Z}_l) := \lim_{\leftarrow} H^1_{\et}({\cV}, {\Z}/{l^r{\Z}}).
$$ 
Furthermore, since ${\cV}$ is a smooth variety, $H_1({\cV}, \Z)$, and hence 
also
${\mt}$, is finitely-generated.  Therefore, there is some integer $t$ such that
$\Hom({\mt}, {\Z}/{l^r{\Z}}) \cong {\Z}/{l^t{\Z}}$ for all $r \geq t$, and 
therefore
$$
\lim_{\leftarrow} \Hom({\mt}, {\Z}/{l^r{\Z}}) \cong {\Z}/{l^t{\Z}}.
$$
\noindent
But ${\Z}/{l^t{\Z}} \otimes \Q = 0$, hence tensoring with ${\Q}$ gives
$$
0 \mapr H^1_{\et}({\alb}_{\cV}, {\Z}_l) \otimes {\Q} 
\mapr H^1_{\et}({\cV}, {\Z}_l) \otimes {\Q} \mapr 0.
$$
\noindent
Since $H^1_{\et}({\cV}, {\Ql}) := H^1_{\et}({\cV}, {\Z}_l) \otimes {\Q}$,
this shows that 
$$
H^1_{\et}({\cV}, {\Ql}) \cong H^1_{\et}({\alb}_{\cV}, {\Ql}).
$$

\noindent
Noting now that there is an isomorphism defined over $k$
$$
\alb(\pvv) \cong \alb_{\cV}, \quad \text {{\cite {sl1}}, 
Thm VI.1.1, p. 148}
$$
\noindent
we see that
\begin{eqnarray*}
H^1_{\et}(\pvv, {\Ql}) & \cong & H^1_{\et}({\alb}(\pvv), {\Ql}) \\
 & \cong & H^1_{\et}({\alb}_{\cV}, {\Ql}) \\
 & \cong & H^1_{\et}({\cV}, {\Ql}).
\end{eqnarray*}
\noindent
as desired. \end{Proof}

\medskip
\noindent
This leads us to the main theorem of this section:

\setcounter{thm}{0}
\medskip
\noindent
\begin{thm}
$$
H^1_{\et}({\cS}/{\bk}, {\Ql}) 
\cong  H^1_{\et}({\cE}/{\bk}, {\Ql}).
$$
\noindent
as $\gal({\bk}/k)$-modules.
\end{thm}
\setcounter{thm}{3}

\par
\noindent
\begin{Proof}  This follows from Theorem {\ref {T3}}, together with 
Proposition {\ref {PVH}}. \end{Proof}

\medskip
\noindent
The following Corollary will be especially useful in later simplifying 
the equation for ${\Ap}$.
 
\medskip
\noindent
\begin{cor}
${\ap}({\cE}) = {\ap}({\cS}).$
\label{C1}
\end{cor}

\par
\noindent
\begin{Proof}.  Since Theorem {\ref {T8}} gives 
$$
H^1_{\et}({\cE}, {\Ql}) \cong H^1_{\et}({\cS}, {\Ql}),
$$
\noindent
it is clear that
$$
\trace(\fr | H^1_{\et}({\cE}/{\bk}, {\Ql})) = 
\trace(\fr | H^1_{\et}({\cS}/{\bk}, {\Ql})).
$$
\noindent
Thus, what remains is to show, for any variety ${\cV}/k$, that
$$
\trace(\fr | H^1_{\et}({\cV}/{\bk}, {\Ql})) = 
\trace(\fr | H^1_{\et}({\tV}/{\bFp}, {\Ql})),
$$
\noindent
and this is a consequence of the proper and smooth base change theorems 
(see {\cite {sga}}). \end{Proof}

\bigskip

\section{A Shioda-Tate Formula}

\medskip
\par
\noindent
In this section, we prove that $\ns{\cE} \otimes \, \Q$ is generated 
by ${\cT}$ and
by the \mbox{$k$-rational} sections.  In the case of elliptic surfaces,
this is the main result of the {\bf Shioda-Tate Formula}, which we 
recall here:

\begin{thm}[Shioda-Tate Formula]{\em ({\cite {ts2}, Theorem 1.1})}

\par
\noindent
Let \mbox{${\mathfrak E} \mapr {\mathfrak C}$} be an elliptic surface defined
over $k$, with \mbox{$k$-rational} section ${\sigma}_0$.
\par
\noindent
Embed ${\mathfrak E}(\mathfrak C/{\bk})$ into 
$\ns{\mathfrak E}$ by sending a section $\sigma$ to the divisor
\mbox{$\sigma(\mathfrak C) - \sigma_0(\mathfrak C)$}.  Then there is a
decomposition of Gal(${\bk}/k$)-modules,
$$
\ns{\mathfrak E} \otimes \Q \cong 
({\mathfrak E}(\mathfrak C/{\bk}) \otimes \Q) \oplus {\cal J},
$$
\noindent
where ${\cal J}$ is the subspace of 
\mbox{$\ns{\mathfrak E} \otimes \Q$} generated by the image of the zero
section, and by all components of all fibers.
\label{TSTF}
\end{thm}

\medskip
\par
\noindent
In order to prove a similar formula for elliptic $n$-folds, we follow 
an argument similar to Shioda's proof for Elliptic Surfaces {\cite {sh1}}, 
and take a 
closer look at the map ${\psi}$ (equation \ref {B9}, restriction to
the generic fiber) at the level of geometric points.

\par
\noindent
Restricting a line bundle on ${\cE}$ to the generic fiber
$E$ defines a homomorphism
\begin{equation}
{\pic}({\cE}) \mapr {\pic}({E/{\hK}})
\label{B4}
\end{equation}

\noindent
which associates with every divisor class cl($D$) on ${\cE}$ the divisor 
\mbox{$D|_{E} = D.E$} on the generic fiber $E$.
Then, using the given rational point 
\mbox{$O \in E(K)$}, adjust the image by sending 
cl($D$) to cl($D'$), where \mbox{$D' := D.E - (D.E)O$}; 
the divisor ${D'}$ is thus a degree zero divisor on $E$, and the 
homomorphism becomes
\begin{equation}
{\phi}: {\pic}({\cE}) \mapr {\pv}(E/{\hK}) \cong E({\hK}).
\label{B5}
\end{equation}

\noindent
\begin{remark}
By construction, the kernel of ${\phi}$ contains the zero-section and any 
fibral divisor in ${\pic}({\cE})$.
\label{R5}
\end{remark}

\par
\noindent
\begin{lemma}
Let ${\kT}$ be the subgroup of ${\pic}({\cE})$ generated by the irreducible
components of the fibral divisors, and by the zero-section $(O)$.  Then
\begin{equation}
0 \mapr {\kT} \maprlim{\eta} {\pic}({\cE}) \maprlim{\phi} E({\hK}) \mapr 0
\label{B7}
\end{equation}

\noindent
is a short exact sequence of Abelian groups.
\label{L1}
\end{lemma}

\noindent
\begin{Proof}
\noindent
\begin{enumerate}
\item
Notice that the morphism ${\phi}$ is surjective: given any 
\mbox{${\hK}$-rational} divisor
$C$ on $E$, taking the schematic closure of its irreducible
components gives a divisor ${\bar C}$ on ${\cE}$ such that
\mbox{${\bar C}.E = C$}.

\item
By Remark (\ref {R5}), \mbox{${\kT} \subset$ ker$(\phi)$}.  
To show that \mbox{${\kT} = ker(\phi)$}, 
consider \mbox{${\Upsilon} \in$ ker$(\phi)$}, 
i.e., \mbox{${\Upsilon} = cl(D)$}, 
where \mbox{$D|_E \sim 0$ on $E$}.
But then \mbox{$D|_E = div(h)$}, where
$$ 
h \in {\hK}(E) = {\bk}({\cS})(E) = {\bk}({\cE}),
$$
\noindent
and hence there exists \mbox{$H \in {\bk}({\cE})$} such that 
\mbox{$(H)|_E = (h)$}.

\par
\noindent
If \mbox{$D' := D - (H)$}, then $D'$ must be in some fiber, i.e.,
\mbox{$D' \in {\kT}$}, and therefore 
$$
{\Upsilon} = cl(D) = cl(D') \in {\kT}.
$$
\end{enumerate} \end{Proof}

\begin{thm}[A Shioda-Tate Formula for Elliptic $N$-folds]

\medskip
\noindent
Embed ${\cE}({\cS}/{\bk})$ into 
$\ns{{\cE}}$ by sending a section $\sigma$ to the divisor
\mbox{${\overline {\sigma({\cS})}} - \sigma_0({\cS})$} (where
\mbox{${\overline {\sigma({\cS})}}$} is the schematic closure
of $\sigma({\cS})$ in ${\cE}$).  

\par
\noindent
Then there is a
decomposition of Gal(${\bk}/k$)-modules,
$$
\ns{{\cE}} \otimes \Q \cong 
({\cE}({\cS}/{\bk}) \otimes \Q) \oplus {\cT}.
$$
\label{T4}
\end{thm}

\par
\noindent
\begin{Proof}
Comparing the short exact sequences (\ref{B6}) and
(\ref{B7}), we see that ${\psi}$ maps ${\pic}({\cE})$ surjectively 
onto $E({\hK})$, while at the same time sending ${\pv}({\cE})$ to 
\mbox{$B({\bk}) = 0$}.  This implies that 
\begin{equation}
\ns{{\cE}} := {\pic}({\cE})/{\pv}({\cE}) \twoheadrightarrow E({\hK}),
\label{B8}
\end{equation}
\noindent
with kernel ${\cT}'$, the image of ${\kT}$ in $\ns{{\cE}}$.
Thus, we have an exact sequence
$$
0 \mapr {\cT}' \mapr \ns{{\cE}} \mapr  E({\hK}) \mapr 0.
$$
\noindent
Since the action of Galois sends vertical divisors to vertical divisors, and 
horizontal to horizontal, this sequence
splits as a Galois module after tensoring with ${\Q}$;
noting that \mbox{${\cT}' \otimes \Q = {\cT}$} then gives 
the desired formula. \end{Proof}

\medskip
\noindent
\begin{cor}
$$
\rank\ns{{\cE}/k} = 1 + \rank{\cE}({\cS}/k) + \rank\ns{{\cS}/k} 
+ \rank{\cF}^{\gal({\bk}/k)},
$$
where ${\cF}$ is the vector space generated by the non-identity 
geometrically irreducible components of the fibral divisors.
\label{C2}
\end{cor}

\medskip
\noindent
\begin{Proof}  Take \mbox{$\gal({\bk}/k)$-invariants} of the Shioda-Tate 
formula for elliptic \mbox{$n$-folds}.  This gives
$$
\rank\ns{{\cE}/k} = \rank{\cE}({\cS}/k) + \rank{\cT}^{\gal({\bk}/k)}.
$$
\noindent
But recall (see Definition {\ref {D5}} ff) that ${\cT}$ is generated by $(0)$, 
${\pi}^*(\ns{\cS})$, and ${\cF}$.

\par
\noindent
Therefore,
$$
\rank{\cT}^{\gal({\bk}/k)} = 1 + \rank\ns{{\cS}/k} 
+ \rank{\cF}^{\gal({\bk}/k)}
$$
\noindent
as required. \end{Proof}

%% file: artsingfib.txt

\section{The Singular Fibers}

The main goal of this section is to prove Theorem {\ref {T5}} below,
which establishes a geometric interpretation for 
the action of Frobenius on the singular fibers.  Our main tools will be 
Tate's Algorithm for determining the singularity-type of a given fiber, and
an effective version of the Geometric Tchebotarev Density Theorem, and this
requires that we now restrict to the case of an {\bf elliptic threefold}
${\cE}/k$.

\begin{thm}  Let ${\cE}/k$ be an elliptic threefold, with notation as 
before.  Then 
$$
{\sum_{x \in {\tD}({\Fp})} (m_x - 1)} = 
{\qp}{\trace({\fr}|{\tF})} + O(\sqrt{\qp}),
$$
where ${\tF}$ is the vector space generated by all non-identity components
of ${\pi}^{-1}({\tD})$ (the identity component is the component intersecting
$(0)$).
\label{T5}
\end{thm}

\par
\noindent
\begin{Proof}

\medskip
\noindent
\begin{enumerate}

\item
Recall (Definition {\ref {DN6}}) that if  
\mbox{${\Delta} = {\Delta}_1 + \cdots + {\Delta}_r$}
is the irreducible decomposition of the discriminant ${\Delta}$
of \mbox{$\pi:{\cE} \mapr {\cS}$},  
then, for all 
\mbox{${\gp} \notin B$}, the discriminant locus of 
\mbox{${\tilde {\pi}}:{\tE} \mapr {\tS}$} is given by
${\tD}$, the reduction of ${\Delta}$ mod(${\gp}$), which has irreducible
decomposition
$$
{\tD} = {\tD}_1 + \cdots + {\tD}_r.
$$ 
\par
\noindent
From this it follows that, for all 
${\gp} \notin B$ (where we recall, from section 2.5, that $B$ is a finite
set of ``bad'' primes), genus ${\tD}_i$ = genus ${\Delta}_i$ for all $i$.  
 
\medskip

\item
Now assume that, for all \mbox{$1 \leq j \leq r$}:
$$
{\sum_{x \in {\tD}_j({\Fp})} (m_x - 1)} =
{\qp}{\trace({\fr}|{\tF}_j)} + O(\sqrt{\qp}),
$$
where ${\tF}_j$ is the vector space generated by the
non-identity fibral divisors over ${\tD}_j$.  Then
\begin{eqnarray*}
{\sum_{x \in {\tD}({\Fp})} (m_x - 1)} & = &
{\sum_{j=1}^r}{\sum_{x \in {\tD}_j({\Fp})} (m_x - 1)}, \\
 & = & {\sum_{j=1}^r} {\qp}{\trace({\fr}|{\tF}_j)} + O(\sqrt{\qp}), \\
 & = & {\qp}{\trace({\fr}|{\tF})} + O(\sqrt{\qp}). \\
\end{eqnarray*}
\noindent
Thus, it suffices to prove the theorem with ${\Delta}$
replaced by one of its irreducible components.

\medskip

\item
The singular points on ${\Delta}$ are determined by the simultaneous 
vanishing of a finite set of polynomials.  But, for a given polynomial
$f(D) \in k[{\Delta}]$, there are only finitely many $x \in {\Delta}$
such that $f(x) \neq 0$, but $f(x) \equiv 0 \mod ({\gp})$ for 
all primes ${\gp}$.  Thus, there can only be finitely-many nonsingular
$x \in {\Delta}$ which are singular points on ${\tD}$ for all reductions
mod(${\gp}$).

\medskip
\noindent
Let $G$ be the set containing all singular points on ${\Delta}$ and
all points which are singular mod(${\gp}$) for all ${\gp}$.

\medskip
\noindent
At the singular points of $s \in {\tD}$, we note that the number
of components lying over ${\tE}_s$ is bounded by the number of
components lying over ${\cE}_s$ (this follows from Section 2.5, 
Statement 4), hence is bounded by a 
constant independent of ${\gp}$. 
Furthermore, the following argument will show that the number of
singular points on ${\tD}$ is bounded independently of ${\gp}$, and
therefore singular points and the fibers over them contribute only to
the error term in Theorem {\ref {T5}}.

\medskip
\noindent
Denote by
\mbox{${\eta}:{\hat {\Delta}} \mapr {\Delta}$} the
desingularization of ${\Delta}$.  If ${\Delta}_{ns}$ is
the set of nonsingular points on ${\Delta}$, then
$$
{\eta}^{-1}({\Delta}_{ns}) \cong {\Delta}_{ns},
$$
\noindent
while, for the finitely-many points $\{s_g\}_1^R \in G$,
$\#{\eta}^{-1}(s_g) = S_g$, a finite constant.

\medskip
\noindent
For almost all ${\gp}$, 
\mbox{${\hat {\tD}} = {\tilde {\hat {\Delta}}}$}, and
\mbox{$\#{\eta}^{-1}(s_g) = \#{\eta}^{-1}({\tilde {s}}_g)$}.
Thus, letting \\
\noindent
\mbox{$S := \max_g \{S_g\}$}, we see that for all
primes ${\gp} \notin B$, 
$$
|\#{\hat {\tD}}({\Fp}) - \#{\tD}({\Fp})| \leq R{\cdot}S.
$$

\medskip

\item
The next step is to show that at almost every point 
\mbox{$x \in {\tD}$}, the singular fiber 
\mbox{${\pi}^{-1}(x)$} is of Kodaira type, and that the singularity type is
locally constant.  This can be accomplished by considering
the localization of ${\tE}$ at ${\tD}$, and applying Tate's Algorithm.

\medskip
\noindent
Let $\tOS$ be the structure sheaf of ${\tS}$, and
${\tOSD}$ the local ring of ${\tD}$ on ${\tS}$.  Then $R := {\tOSD}$ is
a DVR, with residue field \mbox{$F := {{\Fp}({\tD})} = {\Fp}(D)$}, and
prime ${\mathfrak m}$;
let 
$$
{\tE}_{\tD} := {\tE} \times_{\tOS} {\tOSD}
$$
\noindent
be the localization
of ${\tE}$ at ${\tD}$; then a Weierstrass Equation for
${\tE}_{\tD}$ is given by:
\begin{equation}
Y^2 + a_1XY + a_3Y =
X^3 + a_2X^2 + a_4X +a_6, \quad a_i \in R = {\Fp}[D].
\label{S4}
\end{equation}

\par
\noindent
\begin{remark}
Notice that ${\tE}_{\tD}$ is an elliptic curve defined over
a DVR $R$ with {\em non-perfect} residue field $F$, whereas Tate's 
Algorithm is for elliptic curves over DVR's with {\em perfect} residue
fields.  However, the proof works verbatim for any residue field with
the property that all extensions of degree 2 or 3 are separable.  
(For the proof of Tate's Algorithm, see Tate's original paper {\cite {jt2}},
or Silverman {\cite {js} IV.9}).  Thus, if $B$ is expanded to 
include all primes \mbox{$\gp$} such that 
\mbox{$2 | {\qp}$ or $3 | {\qp}$}, then Tate's Algorithm can also be 
applied here.
\label{R6}
\end{remark}

\medskip

\item
Localizing ${\cE}$ over ${\Delta}$ gives an elliptic curve
defined over the DVR ${\cOSD}$, with perfect residue field 
$k({\Delta})$ - hence a straightforward application of Tate's Algorithm
shows that the Kodaira singularity type of ${\cE}$ is locally constant
over ${\Delta}$.

\medskip
\noindent
Then, for all but finitely-many ${\gp}$,

\begin{enumerate}
\item
The Kodaira type of ${\cE}$ over ${\Delta}$ is the same as the
Kodaira type of ${\tE}$ over ${\Delta}$, call it ${\mathfrak D}$.
\medskip

\item
Let
$$
{\Delta}_{KF} := \{ x \in {\Delta} | x \notin G \: 
\text {and ${\cE}_x$ has Kodaira type ${\mathfrak D}$} \}.
$$  
\noindent
Then ${\Delta}_{KF} = {\tD}_{KF}$.

\end{enumerate}

\par
\noindent
The first statement follows because Tate's Algorithm is applied
to the local equations of ${\cE}$ and ${\tE}$ over ${\Delta}$ 
and ${\tD}$ respectively, and for almost all ${\gp}$, they will be
the same equation, with the same zeros.

\medskip
\noindent
The second statement follows from Elimination Theory arguments.

\medskip

\item
The proof now reduces to an examination of each fiber type.

\par
\noindent
Note first that, since we are considering the action of $\fr$ on a
subspace of \mbox{$\ns{{\tE}/{\bFp}}$}, we must look at all fibral
components that are irreducible over ${\mathbf {\bFp}}$ (and not over 
${\bar F}$!).  \mbox{$\trace(\fr)$} picks out from among these the ones that 
are defined over $\Fp$.  These will be called 
\mbox{$\Fp$-rational} components.

\par
\noindent
Following Tate's notation, we now define
$$
a_{i, j} := {\mathfrak m}^{-j}a_i.
$$

\par
\noindent
Consider the divisor ${\tilde {\pi}}^{-1}({\tD})$.
Let 
$$
\begin{array}{lcl}
H & = &  \text {number of irreducible components over ${\bar F}$,} \\
h & = &  \text {number of irreducible components over ${\bar F}$ that are
defined over ${\Fp}$,} \\
M & = &  \text {number of ${\Fp}$-rational components}.
\end{array}
$$
  
\par
\noindent
\begin{itemize}
\item
Type $I_n$, $n > 0$.
\par
\noindent
There are $n$ components of 
\mbox{${\tilde {\pi}}^{-1}({\tD})$} defined over ${\bar F}$, so $H = n$.
According to Tate's Algorithm,
the number $M$ of \mbox{${\Fp}$-rational} components is determined
by the splitting field $F'$ of
$$
P(T) = T^2 + a_1T - a_2.
$$
\noindent
\begin{enumerate}
\item
If \mbox{$F' = F$}, then \mbox{$M = n$};
furthermore, \mbox{$F' = F$} implies that, for all \mbox{$d \in {\tD}$},
$$
P_d(T) = T^2 + a_1(d)T - a_2(d)
$$
\noindent
splits in ${\Fp}$.  Therefore we also have 
\mbox{$m_d = n$} for every \mbox{$d \in {\tD}({\Fp})$}.

\medskip

\item
If \mbox{$F' \neq F$} and \mbox{$F' = L({\tD})$} for 
$L$ some finite extension of $\Fp$, then in fact the $n$ components
of ${\pi}^{-1}({\tD})$ are defined over ${\bFp}$, but they are not all 
defined over ${\Fp}$ because \mbox{$F' \neq F$}.  To determine the number
of ${\Fp}$-rational components, we consider what happens at the fibers
over points \mbox{$d \in {\tD}({\Fp})$}.

\medskip
\noindent
\mbox{$F' = L({\tD})$} implies that, for all \mbox{$d \in {\tD}({\Fp})$},
the splitting field of $P_d(T)$ will be $L \neq {\Fp}$.  Hence, the 
Galois action reflects the geometric components of the polygon.  This
means there will always be two invariant components when $n$ is even,
and one invariant component when $n$ is odd.  Therefore we have:

\par
\noindent
for even $n$: \mbox{$M = m_d = 2$}, for all \mbox{$d \in {\tD}({\Fp})$},

\par
\noindent
for odd $n$: \mbox{$M = m_d = 1$}, for all \mbox{$d \in {\tD}({\Fp})$}.

\item
Otherwise, from the discussion in part (b) it is clear that $h = 2$ for
even $n$, $h = 1$ for odd $n$, and
\begin{eqnarray*}
M & = & {\frac{n+2}{2}} = {\frac{H + h}{2}} \quad \text {for even $n$,} \\
M & = & {\frac{n+1}{2}} = {\frac{H + h}{2}} \quad \text { for odd $n$.} 
\end{eqnarray*}

\par
\noindent
This can be explained as follows:  since \mbox{$F' \neq L({\tD})$},
the irreducible components of ${\pi}^{-1}({\tD})$ are 
NOT defined over ${\bFp}$.  Thus, the ${\bFp}$-components are
exactly the ${\Fp}$-components, and consist of the ${\fr}$ orbits.
Since the action on $I_n$ is reflection, there will be 
\mbox{${\frac{n+2}{2}}$} for even $n$,
and \mbox{${\frac{n+1}{2}}$} for odd $n$.

\medskip
\noindent 
At the level of fibers over points \mbox{$d \in {\tD}({\Fp})$},
$m_d$ alternates between $1$ and $n$ (in the case $n$ odd),
and between $2$ and $n$ (if $n$ even), depending on whether or not
${\Fp}$ is the splitting field of $P_d(T)$. 

\end{enumerate}

\medskip

\item
Type $II$.
There is only one component of 
\mbox{${\tilde {\pi}}^{-1}({\tD})$} defined over ${\bar F}$; since
this includes the identity component ${\Theta}_0$, which is defined
over $F$, we get $M = m_d = 1$ for all \mbox{$d \in {\tD}(\Fp)$}.
\medskip

\item
Type $III$.
There are two components of 
\mbox{${\tilde {\pi}}^{-1}({\tD})$} defined over ${\bar F}$.
One is ${\Theta}_0$, so the other must also remain fixed under the action
of Frobenius. Hence $M = m_d = 2$ for all \mbox{$d \in {\tD}(\Fp)$}.
\medskip

\item
Type $IV$.
There are three components of 
\mbox{${\tilde {\pi}}^{-1}({\tD})$} defined over ${\bar F}$, so $H = 3$.

Let $F'$ be the splitting field of 
$$
P(T) = T^2 + a_{3, 1}T - a_{6, 2}.
$$ 
\medskip

\noindent
\begin{enumerate}
\item
If \mbox{$F' = F$}, then
\mbox{$M = m_d = 3$} for all \mbox{$d \in {\tD}({\Fp})$}.
\medskip

\item
If \mbox{$F' \neq F$} and \mbox{$F' = L(\tD)$}, 
where $L$ is a finite extension of ${\Fp}$,
then \mbox{$M = m_d = 1$} for all \mbox{$d \in {\tD}({\Fp})$}.
\medskip

\item
Otherwise, $m_d$ alternates between $1$ and $3$,
which shows that $h = 1$, and $M = 2 = {\frac{H + h}{2}}$.

\end{enumerate}
\medskip

\item
Type $I_n^*$.
\par
\noindent
There are \mbox{$n + 5$} components of 
\mbox{${\tilde {\pi}}^{-1}({\tD})$} defined over ${\bar F}$.
To determine the number of \mbox{$\Fp$-rational} components, consider the
polynomial
\begin{equation}
Q(T) = T^3 + a_{2, 1}T^2 + a_{4, 2}T + a_{6, 3}.
\label{S5}
\end{equation}

\par
\noindent
\begin{enumerate}
\item
If $Q(T)$ has three distinct roots in ${\bar F}$, then we are in the case
$I_0^*$.  There are 5 components defined over ${\bar F}$
(i.e., $H = 5$.  Let $F'$ be the splitting field of $Q(T)$.

\begin{enumerate}
\item
If $F' = F$, then $M = m_d = 5$ for all $d \in {\tD}(\Fp)$.
\medskip

\item
If $Q(T)$ is reducible in $F$, but $F' \neq F$, then let $P_1(T)$ 
be the irreducible part of $Q(T)$, and let $F''$ be the splitting
field of $P_1(T)$.

\par
\noindent
If $F'' = L({\tD})$, where $L$ is a finite extension of ${\Fp}$,
then $M = M_d = 3$ for all $d \in {\tD}(\Fp)$.

\par
\noindent
If $F'' \neq L({\tD})$, then $M = 4$, and $m_d$ alternates
between 3 and 5.
\medskip

\item
If $Q(T)$ is irreducible in $F$, and 
$F' = L({\tD})$, where $L$ is a finite extension of ${\Fp}$,
then $M = m_d = 2$ for all $d \in {\tD}(\Fp)$.

\item
If $Q(T)$ is irreducible in $F$, and 
$F' \neq L({\tD})$, then $m_d$ alternates
between 2, 3, and 5, (so $h = 2$), and $M = 3 = {\frac{H + h}{2}}$.

\end{enumerate}
\medskip

\item
If $Q(T)$ has one simple root and one double root in ${\bar F}$, then
we have type $I_n^*, n \geq 1$.  In this case, $H = n + 5$, and
the number of ${\Fp}$-rational components
(\mbox{$n+3$} or \mbox{$n + 5$}) is determined by the splitting of a
quadratic polynomial $P(T)$.  Let $C$ be its splitting field.

\begin{enumerate}
\item
If $C = F$, then $M = m_d = n + 5$ for all $d \in {\tD}(\Fp)$.
\medskip

\item
If $C = L({\tD})$, where $L$ is a finite extension of ${\Fp}$,
then $M = m_d = n + 3$ for all $d \in {\tD}(\Fp)$.
\medskip

\item
If $C \neq L({\tD})$, then $m_d$ alternates
between $n + 3$, and $n + 5$ (so $h = n + 3$), and 
$M = n + 4 = {\frac{H + h}{2}}$.

\end{enumerate}

\end{enumerate}

\medskip

\item
Type $III^*$.
\par
\noindent
There are eight components of 
\mbox{${\tilde {\pi}}^{-1}({\tD})$} defined over ${\bar F}$, and
they are all fixed by the action of Frobenius.  Hence, 
\mbox{$M = m_d = 8$} for all $d \in {\tD}(\Fp)$.
\medskip

\item
Type $II^*$.
\par
\noindent
There are nine components of 
\mbox{${\tilde {\pi}}^{-1}({\tD})$} defined over ${\bar F}$, and
they are all fixed by the action of Frobenius.  Hence, 
\mbox{$M = m_d = 9$} for all $d \in {\tD}(\Fp)$.
\medskip

\item
Type $IV^*$.
\par
\noindent
There are seven components defined over ${\bar F}$, therefore $H = 7$.

\par
\noindent
Let $F'$ be the splitting field of
$$
P(T) = T^2 + a_{3, 2}T - a_{6, 4}.
$$

\begin{enumerate}
\item
If $F' = F$ then \mbox{$M = m_d = 7$} for all $d \in {\tD}(\Fp)$.
\medskip

\item
If $F' = L({\tD})$, where $L$ is a finite extension of ${\Fp}$,
then $M = m_d = 3$ for all $d \in {\tD}(\Fp)$.
\medskip

\item
If $F' \neq L({\tD})$, then $m_d$ alternates
between 3 and 7 (so $h = 3$), while $M = 5 = {\frac{H + h}{2}}$.

\end{enumerate}

\end{itemize}

\medskip
\item
The proof of the theorem now reduces to analyzing two cases:
\begin{enumerate}
\item
The case $M = m_d$ for all  $d \in {\tD}(\Fp)$:

\par 
\noindent
We have
\begin{eqnarray}
\sum_{d \in {\tD}(\Fp)} (m_d - 1) & = &
\sum_{d \in {\tD}(\Fp)} (M - 1) \nonumber \\
 & = & \#{\tD}(\Fp)(M - 1)   \nonumber \\
 & = & \#{\tD}(\Fp){\trace({\fr}|{\tF})} 
\label{SF1}
\end{eqnarray}

\medskip
\item
The case $M \neq m_d$:
This requires a Chebotarev argument.

\medskip

We use the following effective version of the geometric Chebotarev Density
Theorem:

\begin{thm}[Murty, Scherk {\cite {ms}}]
If \mbox{$X \mapr Y$} is a geometric covering of curves over ${\Fp}$
(i.e., $Y$ is defined over ${\Fp}$, and ${\Fp}$ is the algebraic
closure of itself in ${\Fp}(X)$), then
\begin{equation}
\left|\psi_C(r) - {\frac{|C|}{|G|}}\psi(r)\right| \leq
2g_X{\frac{|C|}{|G|}}{\qp}^{\frac{r}{2}} + |D|;
\end{equation}
\noindent
where, for \mbox{$y \in Y$} unramified,
$C$ is a conjugacy class in $G := Gal(X/Y)$,
$\sigma_y$ is a Frobenius conjugacy class of $y$,
and 
$$
\begin{array}{lcl}
Y_r       & = & Y \times_{\Fp} {\F}_{{\gp}^r}, \\
{\bar Y}  & = & Y \times_{\Fp} {\bFp}, \\
\psi_C(r) & = & \#\{y \in Y_r |\text {$y$ unramified, $\sigma_y \in C$} \}, \\
\psi(r)   & = & \#\{y \in Y_r |\text {$y$ unramified}\}, \\
D         & = & \text {the set of ramified points in ${\bar Y}$.} \\
\end{array}
$$
\label{TMS}
\end{thm}

\medskip

\item
To apply this version of Chebotarev, we note that in our case,
$X$ is the curve defined by the polynomial $P(T)$, $Y = {\tD}$, and,
since we are considering only points in ${\tD}(\Fp)$, $r = 1$.
Thus, we get the inequality:
\begin{eqnarray*}
|\psi_C(1) - {\frac{|C|}{|G|}}\psi(1)| & \leq & 
2g_X{\frac{|C|}{|G|}}{\sqrt{\qp}} + |D|, \\
& \leq & 2g_X{\sqrt{\qp}} + |D|. \\
\end{eqnarray*}
\noindent
Since, furthermore, by the Riemann-Hurwitz formula
$$
2g_X - 2 = (2g_Y - 2)|G| + |D|,
$$
\noindent
the inequality becomes
\begin{eqnarray*}
|\psi_C(1) - {\frac{|C|}{|G|}}\psi(1)| & \leq & 
2g_Y{\sqrt{\qp}} 
+ (2 - 2|G| + |D|){\sqrt{\qp}} + |D|; \\
& \leq & K{\sqrt{\qp}} + |D|. \\
\end{eqnarray*}

where $K$ is a constant depending only on $|G|$, $g_Y$ and $|D|$, and hence is
bounded independently of ${\gp}$.
This gives
\begin{equation}
\psi_C(1) = {\frac{|C|}{|G|}}\psi(1) + O({\sqrt{\qp}}).
\end{equation}
\noindent
Noting that 
\mbox{$\psi(1) = \#{\tD}(\Fp) - |D|$},
gives
$$
\psi_C(1) = {\frac{|C|}{|G|}}\#{\tD}(\Fp) + O({\sqrt{\qp}}).
$$
\medskip
\noindent

Next, we want to observe what happens for deg(P) = 2, and deg(P) = 3.

\medskip
\noindent
\begin{enumerate}
\item
deg(P) = 2.

\medskip
\noindent
In this case, $|G| = 2$; one of the two conjugacy classes consists of
$id$, and the other is a transposition $t$.

\par
\noindent
If, for a given $d \in {\tD}$, $\sigma_d = id$, then $m_d = H$.  If
$\sigma_d = t$, then $m_d = h$.
Thus, the sum of components of fibers breaks up into:
\begin{eqnarray*}
\sum_{d \in {\tD}(\Fp)} m_d & = & \sum_{d: \sigma_d = id} m_d
+ \sum_{d: \sigma_d = t} m_d \\
 & = & \sum_{d: \sigma_d = id} H
+ \sum_{d: \sigma_d = t} h \\
 & = & H|\psi_{id}| + h|\psi_{t}|\\
 & = & H{\frac{1}{2}}\#{\tD}(\Fp) + h{\frac{1}{2}}\#{\tD}(\Fp) 
+ O({\sqrt{\qp}}) \\
 & = & {\frac{H + h}{2}}\#{\tD}(\Fp) + O({\sqrt{\qp}}). 
\end{eqnarray*}
\medskip

\item
deg(P) = 3.

\medskip
\noindent
In this case, $|G| = 3$ or $|G| = 6$. 

\noindent
\begin{enumerate}
\item
If $|G| = 3$, there are three conjugacy classes corresponding to the three 
elements of $G$: $id$, $v_1$, and $v_2$.

\par
\noindent
If, for a given $d \in {\tD}$, $\sigma_d = id$, then $m_d = 5$; if
$\sigma_d = v_i$, $i = 1, 2$, then $m_d = 2$.
Thus, the sum of components of fibers breaks up into:
\begin{eqnarray*}
\sum_{d \in {\tD}(\Fp)} m_d & = & \sum_{d: \sigma_d = id} m_d
+ \sum_{d: \sigma_d = v_1} m_d + \sum_{d: \sigma_d = v_2} m_d\\
 & = & \sum_{d: \sigma_d = id} 5
+ \sum_{d: \sigma_d = v_1} 2 + \sum_{d: \sigma_d = v_2} 2\\
 & = & 5|\psi_{id}| + 2|\psi_{v_1}| + 2|\psi_{v_2}|\\
 & = & 5({\frac{1}{3}})(\#{\tD}(\Fp)) + 2({\frac{1}{3}})(\#{\tD}(\Fp)) \\
 &   & + 2({\frac{1}{3}})(\#{\tD}(\Fp)) + O({\sqrt{\qp}}) \\
 & = & 3\#{\tD}(\Fp) + O({\sqrt{\qp}}). 
\end{eqnarray*}

\medskip
\item
If $|G| = 6$, then there are three conjugacy classes;
one conjugacy class consists of
$id$; a second one is a transposition $t$, and this has 3 elements; the third
conjugacy class is the even permutations $v$, of which there are two.

\par
\noindent
If, for a given $d \in {\tD}$, $\sigma_d = id$, then $m_d = 5$.  If
$\sigma_d = t$, then $m_d = 3$, and if $\sigma_d = v$, then $m_d = 2$.
Thus, the sum of components of fibers breaks up into:
\begin{eqnarray*}
\sum_{d \in {\tD}(\Fp)} m_d & = & \sum_{d: \sigma_d = id} m_d
+ \sum_{d: \sigma_d = t} m_d + \sum_{d: \sigma_d = v} m_d\\
 & = & \sum_{d: \sigma_d = id} 5
+ \sum_{d: \sigma_d = t} 3 + \sum_{d: \sigma_d = v} 2\\
 & = & 5|\psi_{id}| + 3|\psi_{t}| + 2|\psi_{v}|\\
 & = & 5({\frac{1}{6}})(\#{\tD}(\Fp)) + 3({\frac{3}{6}})(\#{\tD}(\Fp))\\
 &   & + 2({\frac{2}{6}})(\#{\tD}(\Fp)) + O({\sqrt{\qp}}) \\
 & = & 3\#{\tD}(\Fp) + O({\sqrt{\qp}}). 
\end{eqnarray*}
\end{enumerate}
\end{enumerate}

\medskip
\item
By Weil's Estimate, we have:
$$
\#{\tD}({\Fp}) = 1 + {\ap}({\Delta}) + {\qp},
$$
\noindent
where, if $g_{\Delta}$ is the genus of ${\Delta}$, 
$$
|{\ap}({\Delta})| \leq {\sqrt{\qp}}g_{\Delta}.
$$
\noindent
Furthermore, ${\trace({\fr}|{\tF})} < M$ is bounded independently
of ${\gp}$, because for almost all ${\gp}$, $M$ is given by the
number of irreducible components of the singular fiber over ${\Delta}$.
Therefore,
\begin{eqnarray*}
\#{\tD}({\Fp}){\trace({\fr}|{\tF})} & = & 
(1 + {\ap}({\Delta}) + {\qp}){\trace({\fr}|{\tF})} \\
 & = & {\qp}{\trace({\fr}|{\tF})} \\
 &   & + (1 + {\ap}({\Delta})){\trace({\fr}|{\tF})} \\
 & = & {\qp}{\trace({\fr}|{\tF})} + O({\sqrt{\qp}}).
\end{eqnarray*}

Thus, combining this with the result in equation (\ref {SF1}) gives
$$ 
\sum_{d \in {\tD}(\Fp)} (m_d - 1) = 
{\qp}{\trace({\fr}|{\tF})} + O(\sqrt{\qp}). 
$$

\end{enumerate}
\end{enumerate}
\end{Proof}

%% file: artlser.txt

\section{The Main Theorem}

\par
\noindent
We now have all the information necessary to prove Theorem {\ref {MThm}}.
What remains is to run it through the $L$-series machinery, and apply 
Tate's Conjecture.

\medskip
\noindent
Recall that, for a smooth variety ${\cV}/k$, 
the {\bf Hasse-Weil ${\mathbf L}$-series} attached to 
$H^i_{\et}({\cV}/{\bk})$, denoted $L_2({\cV}, s)$, is given by
$$
L_2({\cV}, s) := \prod_{\gp} P_{{\gp}, 2}({\tV}, {\qp}^{-s}).
$$
\noindent
where
$$
P_{{\gp}, i} = \det\left(1 - {\fr}t | H^i_{\et}({\cV}/{\bk}; {\Ql})\right).
$$

\medskip
\begin{remark}
To be precise, since in this paper we are working over all primes
${\gp} \in U$ (i.e., we are excluding the primes in $B$),  
\begin{eqnarray*}
L_2({\cE}, s) & \approx & \prod_{\gp} P_{{\gp}, 2}, 
({\tE}, {\qp}^{-s})\quad \text {and} \\
L_2({\cS}, s) & \approx & \prod_{\gp} P_{{\gp}, 2}({\tS}, {\qp}^{-s}),
\end{eqnarray*}
\noindent
where the symbol $\approx$ is used to indicate that the two sides
agree up to finitely many Euler factors.  This, however, has no
effect on the residue computation.
\end{remark}

\medskip
\noindent
\begin{conj}[Tate's Conjecture] (\cite {jt1}, Conjecture 2).
\medskip
\noindent
{\it Let ${\cV}$ be a smooth projective variety defined over $k$, and
let $L_2({\cV}, s)$ be the Hasse-Weil \mbox{$L$-function} attached to
$H^2_{\et}({\cV}/{\bk})$.  Then $L_2({\cV}, s)$ has a meromorphic continuation
to ${\C}$, and has a pole at $s = 2$ of order:
$$
-\ord_{s = 2} L_2({\cV}, s) = \rank \ns{{\cE}/k}.
$$ }
\label{TC}
\end{conj}

\medskip
\noindent
Finally, we are ready to prove the main theorem:

\setcounter{section}{1}
\setcounter{thm}{0}
\medskip
\noindent 
\begin{thm}
Let ${\cE} \rightarrow {\cS}$ be a non-split
elliptic threefold defined over a number field $k$.  Then Tate's Conjecture
for ${\cE}/k$ and ${\cS}/k$ implies
$$ \res_{s=1} \sum_{\gp}{-{\Ap}{\frac{\log{\gp}}{{\gp}^s}}}
= \rank{\cE}({\cS}/k).$$
\end{thm}
\setcounter{section}{6}

\clearpage

\medskip
\par
\noindent
\begin{Proof}

\begin{enumerate}
\item
Look at ${\tE}$ as a fibration of curves, and use the 
Lefschetz Fixed-Point Theorem to count its rational points fiber by fiber:
\begin{eqnarray}
\#{\tE}({\Fp}) & = & {\sum_{x \in {\tS}({\Fp})} \#{\tE}_x({\Fp})} \nonumber \\
& = & {\sum_{x \in {\tS}({\Fp})}{(1-{\ap}({\tE}_x)+{\qp}+(m_x-1){\qp}})}
\nonumber \\
& = & (1+{\qp})\#{\tS}({\Fp}) - {{\qp}^2}{\Ap} + {\sum_{x
\in{\tS}({\Fp})} (m_x - 1){\qp}}.
\label{E1}
\end{eqnarray}

\noindent
Since $m_x = 1$ for all non-singular fibers
${\tE}_x$ with \mbox{$x \in {\tS}(\Fp)$},
$$
{\sum_{x \in{\tS}({\Fp})} (m_x - 1){\qp}} =
{\sum_{x \in{\tD}({\Fp})} (m_x - 1){\qp}}.
$$

\noindent
Therefore, combining equation (\ref {E1}) with the result of 
Theorem {\ref {T5}} gives
\begin{eqnarray}
\#{\tE}({\Fp}) & = & (1+{\qp})\#{\tS}({\Fp}) - {{\qp}^2}{\Ap} \nonumber \\
 & & + \, {\qp}{\trace({\fr}|{\tF}){\qp}} + O(\sqrt{{\qp}^3}).
\label{E7}
\end{eqnarray}

\medskip

\item
Now use the Lefschetz Fixed-Point Theorem again, this time considering
${\cE}$ as a threefold, and ${\cS}$ as a surface:
\begin{eqnarray}
\#{\tE}({\Fp}) & = & 1-{\ap}({\tE})+{\bp}({\tE})-{\cp}({\tE})+{\qp}{\bp}({\tE})
-{\qp}^2{\ap}({\tE})+{\qp}^3.
\label{E2} \\
\#{\tS}({\Fp}) & = & 1-{\ap}({\tS})+{\bp}({\tS})-{\qp}{\ap}({\tS})+{\qp}^2,
\label{E3}
\end{eqnarray}
\noindent
where we note that 
\begin{eqnarray*}
\trace({\fr}|H^5_{\et}({\cE}, \Ql)) & = & {\qp}^2{\ap}({\tE}), \\
\trace({\fr}|H^4_{\et}({\cE}, \Ql)) & = & {\qp}{\bp}({\tE}), \\
\trace({\fr}|H^3_{\et}({\cS}, \Ql)) & = & {\qp}{\ap}({\tS}),
\end{eqnarray*}
are given by duality
\begin{eqnarray*}
H^1_{\et}({\cE}, \Ql) & \cong & H^5_{\et}({\cE}, \Ql){\hat { }} \: , \\
H^2_{\et}({\cE}, \Ql) & \cong & H^4_{\et}({\cE}, \Ql){\hat { }} \: , \\
H^1_{\et}({\cS}, \Ql) & \cong & H^3_{\et}({\cS}, \Ql){\hat { }} \: .
\end{eqnarray*} 

\medskip
\noindent
Since ${\ap}({\cS})={\ap}({\cE})$ by Corollary (\ref{C1}), 
equation (\ref {E2}) implies:
\begin{equation}
\#{\tE}({\Fp}) = 1-{\ap}({\tS})+{\bp}({\tE})-{\cp}({\tE})+{\qp}{\bp}({\tE})
-{{\qp}^2}{\ap}({\tS})+{\qp}^3. 
\label{E5}         
\end{equation}
                                      
\par
\noindent
A different expression for the number of rational points on ${\tE}$ is
obtained by combining equations (\ref {E7}) and (\ref {E3}):
\begin{eqnarray}
\#{\tE}({\Fp}) & = & 1 + {\qp} + {\qp}^2 + {\qp}^3 -{\ap}({\tS}) - 2{\qp}{\ap}
({\tS}) - {{\qp}^2}{\ap}({\tS}) + {\bp}({\tS}) + {\qp}{\bp}({\tS}) \nonumber \\
& & - {{\qp}^2}{\Ap} +{\qp}{\trace({\fr}|{\tF}){\qp}} 
+ O(\sqrt{{\qp}^3}).
\label{E4}
\end{eqnarray}

\medskip

\item
Finally, equating the two expressions for the number of rational points
on ${\tE}$ in equations (\ref {E5}) and (\ref {E4}) gives an expression
for ${\Ap}$:
 
\begin{eqnarray}
{{\qp}^2}{\Ap} & = & {\qp} - 2{\qp}{\ap}({\tS}) + {\bp}({\tS})
+ {\qp}{\bp}({\tS}) + {\cp}({\tE}) - {\bp}({\tE}) - {\qp}{\bp}({\tE}) 
\nonumber \\
& & + {\qp}^2 + {\qp}{\trace({\fr}|{\tF}){\qp}} 
+ O(\sqrt{{\qp}^3}).
\label{E6}
\end{eqnarray}
 
\noindent
By Deligne's Theorem {\cite {pd1}}, we know, for every smooth, 
projective variety ${\cV}$ defined over $k$, that 
\begin{eqnarray*}
|{\ap}({\cV})| & \leq & B_1({\cV}){\sqrt{\qp}}, \\
|{\bp}({\cV})| & \leq & B_2({\cV}){{\qp}}, \\
|{\cp}({\cV})| & \leq & B_3({\cV}){\sqrt{{\qp}^3}},
\end{eqnarray*}
\noindent
where \mbox{$B_i({\cV}) := \dim H^i_{\et}({\cV}/{\bk}, {\Ql})$}
is independent of $\gp$.

\par
\noindent
Thus, we can group all terms of order ${\sqrt{{\qp}^3}}$ or less together,
and obtain:
\begin{equation}
{{\qp}^2}{\Ap}  = {\qp}^2 + {\qp}{\bp}({\tS}) - {\qp}{\bp}({\tE}) 
 + {\trace({\fr}|{\tF}){\qp}^2} + O(\sqrt{{\qp}^3}).
\label{E8}
\end{equation}

\medskip

\item
And now compute residues.

\begin{enumerate}
\item
\begin{equation}
\res_{s=1} \sum_{\gp} {\frac {\log{\qp}}{{\qp}^s}} = 1. \qquad \qquad
\label{E10}
\end{equation}

\medskip

\item
Letting $L({\cF}, s)$ be the Artin $L$-series attached to ${\cF}$, we
have, for \\
\mbox{$\re(s) > \frac{1}{2}$}
\begin{eqnarray*}
\frac{{\rm d}}{{\rm {ds}}} \log{L({\cF}, s)} & = & 
\frac{{\rm d}}{{\rm {ds}}} \sum_{\gp} -\log 
\det\left(1 - {\fr}{\qp}^{-s} | {\cF}\right) \\
 & = & \sum_{\gp} -\trace\left(\fr | {\cF}\right)\frac{\log{\qp}}{{\qp}^s}
+ O(1) .
\end{eqnarray*}
\noindent
Therefore,
\begin{eqnarray}
\res_{s=1} \sum_{\gp} \trace({\fr}|{\tF}){\frac {\log{\qp}}{{\qp}^s}}
& = & -\res_{s = 1}\frac{{\rm d}}{{\rm {ds}}} \log{L({\cF}, s)} \nonumber \\ 
& = & - \ord_{s=1} L({\cF}, s) \nonumber \\
 & = & \rank({\cF}^{\gal({\bk}/k)}).
\label{E11}
\end{eqnarray}
\noindent
This last equality follows from {\cite {rs}}, Proposition 1.5.1. 
\medskip

\item
Furthermore, for $\re(s) > \frac{3}{2}$,
\begin{eqnarray*}
\frac{{\rm d}}{{\rm {ds}}} \log{L_2({\cE}, s)} & = & 
\frac{{\rm d}}{{\rm {ds}}} \sum_{\gp} -\log 
\det\left(1 - {\fr}{\qp}^{-s} | H^2_{\et}({\cE}/{\bk}, \Ql)\right) \\
 & = & \sum_{\gp} -{\bp}({\cE})\frac{\log{\qp}}{{\qp}^s}
+ O(1) .
\end{eqnarray*}
\noindent
Therefore,
\begin{eqnarray}
\res_{s=2} \sum_{\gp} {\bp}({\cE}){\frac {\log{\qp}}{{\qp}^s}}
& = & -\res_{s = 1}\frac{{\rm d}}{{\rm {ds}}} \log{L_2({\cE}, s)} \nonumber\\  
& = & -\ord_{s=2} L_2({\cE}, s) \nonumber \\
 & = & \rank \ns{{\cE}/k} \quad \text {by Tate's Conjecture},
\label{E13}
\end{eqnarray}
\noindent
and similarly
\begin{eqnarray}
\res_{s=2} \sum_{\gp} {\bp}({\cS}){\frac {\log{\qp}}{{\qp}^s}}
& = & -\res_{s = 1}\frac{{\rm d}}{{\rm {ds}}} \log{L_2({\cS}, s)} \nonumber\\  
& = & -\ord_{s=2} L_2({\cS}, s) \nonumber \\
 & = & \rank \ns{{\cS}/k} \quad \text {by Tate's Conjecture}.
\label{E12}
\end{eqnarray}

\end{enumerate}

\medskip

\item
Combining the residue calculations with equation (\ref {E8}) gives
$$
\res_{s=1} \sum_{\gp} -{\Ap}{\frac {\log{\qp}}{{\qp}^s}}
 =  -1 - \rank({\cF}^{\gal({\bk}/k)}) - \rank \ns{{\cS}/k} 
+ \rank \ns{{\cE}/k}.
$$
\noindent
and by the Shioda-Tate formula for elliptic threefolds 
(Theorem {\ref {T4}}), this gives
\begin{equation}
\res_{s=1} \sum_{\gp} -{\Ap}{\frac {\log{\qp}}{{\qp}^s}}
= \rank {\cE}({\cS}/k).
\end{equation}

\end{enumerate}
\end{Proof}

%% file: artnot.txt

\section{Notation}

$$
\begin{array}{ll}
k/{\Q}      & \text {a number field.} \\
{\bk}       & \text {the algebraic closure of $k$.} \\
{\Fp}       & \text {the residue field of a prime $\gp$ of $k$.} \\
{\bFp}       & \text {the algebraic closure of ${\Fp}$.} \\
{\qp}       & \text {the norm of $\gp$, i.e., ${\qp}=\#{\Fp}$.} \\
{\cS}/k     & \text {a smooth, projective surface defined over $k$.} \\
{\cE}/k     & \text {a (non-split) elliptic threefold 
$\pi:{\cE} \rightarrow S$ with section $\sigma: {\cS} \rightarrow {\cE}$;}\\
            & \text {regular, proper and flat over ${\cS}$, and 
defined over $k$.} \\
{\tS}/{\Fp} & \text {the reduction of ${\cS} \pmod{\gp}$.} \\
{\tE}/{\Fp} & \text {the reduction of ${\cE} \pmod{\gp}$.} \\
K           & = k({\cS}) \text {, the function field of ${\cS}/k$.} \\
{\hK}       & = {\bk}({\cS}) \text {, the function field of ${\cS}/{\bk}$.} \\
E/K         & \text {the generic fiber of $\pi:{\cE} \rightarrow {\cS}$, 
a smooth elliptic curve defined over $K$.}\\
({\tau}, B) & \text {the ${\hK}/{\bk}$-trace of $E$.} \\
\Delta      & \text {the discriminant locus of ${\cE} \mapr {\cS}$,} \\
            & \text {with irreducible decomposition 
$\sum_{j=1}^{r}{\Delta}_j$.} \\
{\tD}       & \text {the discriminant locus of the fibration $\pmod{\gp}$,} \\
            & \text {with irreducible decomposition 
$\sum_{j=1}^{r}{\tD}_j$.} \\
M_j         & \text {the number of ${\Fp}$-rational components of 
${\tilde \pi}^{-1}({\tD}_j)$.} \\
{\Theta}_{i,j} & (0 \leq i \leq {M_j - 1}): \text {the ${\Fp}$-rational 
components of ${\tilde \pi}^{-1}({\tD}_j)$.} \\
{\Theta}_{0,j} & \text {the identity component.} \\
m_x         & \text {the number of ${\Fp}$-rational components of the 
fiber ${\tE}_x$, for given $x \in {\tS}({\Fp})$.}
\end{array}
$$